\documentclass[12pt, openright, a4paper, english, brazil]{amsart}
\usepackage[T1]{fontenc}
\usepackage[utf8]{inputenc}
\usepackage{scalefnt}
\usepackage[bottom=2cm,top=2cm,left=2cm,right=2cm]{geometry}
\usepackage{xcolor}
\newcommand\crule[3][black]{\textcolor{#1}{\rule{#2}{#3}}}

\usepackage{amsthm}
\usepackage{amsmath}

\usepackage{amsfonts}
\usepackage{amssymb}
\usepackage{array}
\usepackage{graphicx}
\usepackage{graphics}
\usepackage[all]{xy}
\usepackage{tikz}
\usepackage{subfig}
\usetikzlibrary{matrix}
\usepackage{tikz-cd}
\usepackage[backend=biber, style=alphabetic, sorting=nty]{biblatex}
\usepackage{hyperref}
\usepackage{changepage}
\usepackage{enumitem}
\usepackage[nameinlink]{cleveref}
\usepackage{mathtools}

\newtheorem{definicao}{Definition}[section]
\newtheorem{proposicao}{Proposition}[section]

\newtheorem{teorema}{Theorem}[section]
\newtheorem{observacao}{Remark}[section]

\DeclareUnicodeCharacter{0301}{\'{e}}

\title{Operator $\Delta-aS$ on warped product manifolds}
\author[Barbosa, Souza and Viana]{Ezequiel Barbosa, Mateus Souza and Celso Viana}
\address{Departamento de Matemática, Instituto de Ciências Exatas, Universidade Federal de Minas Gerais}

\begin{document}

\email{ezequiel@mat.ufmg.br}
\email{mateushenriquer@yahoo.com.br}
\email{celso@mat.ufmg.br}

\begin{abstract}
    In this work we studied the stability of the family of operators $L_a=\Delta-aS$, $a\in\mathbb R$, in a warped product of an infinite interval or real line by one compact manifold, where $\Delta$ is the Laplacian and $S$ is the scalar curvature of the resulting manifold. 
\end{abstract}

\maketitle

\section{Introduction}

\addcontentsline{toc}{chapter}{Introduction}

Let $M^n$ be a complete Riemannian manifold of dimension $n\ge 2$, consider the operator $\Delta-q$, where $\Delta$ is the rough Laplacian of $M$ and $q:M\rightarrow\mathbb R$ is a smooth function. D. Fischer-Colbrie and R. Schoen \autocite{fischer1980structure} did a study about this type of operator and concluded that the existence of a positive function $f$ on $M$ satisfying $\Delta f-qf=0$ is equivalent to the condition that the first eigenvalue of $\Delta-q$ under the Dirichlet boundary condition is positive in each limited domain of $M$. 

We say that an operator $L=\Delta-q$ is \textit{stable} if $L$ satisfies
\begin{equation}
\label{eqest1}\int_M-fLf\ge 0
\end{equation}
\noindent for all $f\in C_c^{\infty}(M)$. $L$ is said to be \textit{unstable} when $L$ is not stable. The inequality (\ref{eqest1}) is equivalent to 
\begin{equation}
\label{eqest2}\int_M|\nabla f|^2+qf^2\ge 0
\end{equation}

\noindent for all $f\in C_c^{\infty}(M)$. 
Using arguments of approximation in $H^1$ norm, the space of test functions for (\ref{eqest2}) can be replaced by the space $C_c^{0,1}(M)$ of Lipschitz functions of compact support in $M$. 

Let $D\subset M$ be a bounded domain. According to the theory of elliptical equations, the operator $\Delta-q$ acting on functions with Dirichlet boundary conditions of $D$ has a discrete spectrum $\lambda_1^{(D)}<\lambda_2^{(D)}\le\lambda_3^{(D)}\le\dots$ of eigenvalues. The usual characterization of the first eigenvalue of $\Delta-q$ on $D$ is $$\lambda_1^{(D)}=\inf\left\{\int_D|\nabla f|^2+qf^2; \hspace{2mm} f\in C^{\infty}(M),\hspace{2mm} \text{supp} f\subset D, \int_Df^2=1\right\}.$$

We can conclude that an operator $L=\Delta-q$ is stable if and only if the first eigenvalue of $L$ is positive on each bounded domain under the Dirichlet boundary condition.

Let $S$ be the scalar curvature of $M$. We can consider the family of operators $L_a\coloneqq \Delta-aS$, $a\in\mathbb R$. For some values of $a$, there are interesting geometric properties and results, for example: 

\noindent (i) When $a=0$, we have the usual Laplace operator.

\noindent (ii) When $a=\frac{n-2}{4(n-1)}$, $n\ge 3$, the operator $L_a$ is related to the Yamabe operator $Y$ of $M$ by the relationship $Y=\frac{1}{a}L_a$. In particular, if the first eigenvalue of $L_a$ is negative, there is a metric in $M$ conformal to the original metric of $M$ that has constant positive scalar curvature. If the first eigenvalue is zero, then the metric will have constant nonnegative scalar curvature \autocite{li2023metrics}.

\noindent (iii) When $a=\frac{1}{4}$, the operator $L_a$ appears in Perelman's work on three dimensional Ricci flow with surgery \autocite{perelman2002entropy,perelman2003ricci}.

\noindent (iv) When $a=\frac{1}{2}$, $n\ge 3$ and the first eigenvalue of $L_a$ is negative (resp.
nonpositive), there is an isometric immersion of $M$ in a manifold $N$ of positive (resp. nonnegative) scalar curvature such that $M$ becomes
a two-sided stable minimal hypersurface. More precisely, $N$ is diffeomorphic to $M\times\mathbb S^1$ and $L_a$ is the Jacobi operator, referring to the formula for the second variation of the area of $M$ \autocite{li2023metrics}.

\noindent (v) When $a=1$, the operator $L_a$ becomes the Jacobi operator of the second variation of the area of a minimal hypersurface of a flat space.

\noindent (vi) It follows from the definition that, if the first eigenvalue of $L_a$ is negative (resp. nonpositive) for all $a>0$, the scalar curvature of $M$ is positive (resp. nonnegative). Similarly, if the first eigenvalue of $L_a$ is negative (resp. nonnegative) for all $a<0$, the scalar curvature of $M$ is negative (resp. nonpositive).

The main goal of this work is to study the stability of the operator $L_a$ in certain types of Riemannian manifolds $M$, as well as some properties obtained at $M$ from the results to be addressed. As a convention, we say that a Riemannian manifold $M$ is \textit{$a$-stable} if the operator $L_a$ is stable.

\section{Case \texorpdfstring{$n=2$}{j}}\label{888}

The case $n=2$ is a special case of the general case, as the scalar curvature is twice the sectional curvature, which is simply the Gaussian curvature. Hence, the sectional and scalar curvatures are equivalent objects. 
Therefore, the family of operators becomes $L_a=\Delta-2aK$, where $K$ denotes the Gaussian curvature of $M$.  A special case studied by Kawai in \autocite{kawai1988operator} is when $M$ has nonnegative Gaussian curvature. He proved:

\noindent \textbf{Theorem} [Kawai 1984]
\hypertarget{Kawai}{\textit{Let}} \textit{$M$ be an oriented complete noncompact bidimensional Riemannian manifold of nonpositive Gaussian curvature $K$, where $K$ is not identically zero. Suppose that $a>\frac{1}{8}$, then there is a function $f$ of compact support that satisfies the inequality} $$\int_M|\nabla f|^2+2aKf^2<0.$$

In particular, if $M$ has nonpositive Gaussian curvature and is $a$-stable for some $a>\frac{1}{8}$, then $M$ is flat, that is, its Gaussian curvature is identically zero. Furthermore, under the condition $K\le 0$, is trivial in that $M$ is $a$-stable for all $a\le 0$. With that, if the Gaussian curvature of $M$ is nonpositive, we have an "interval of stability," being it $(-\infty,0]$ and an "interval of instability," being it $(\frac{1}{8},\infty)$. Thus, there is an interval of uncertainty regarding the stability of $L_a$ for a two-dimensional manifold of nonpositive Gaussian curvature, being it $(0,\frac{1}{8}]$. In fact, the flat plane is $a$-stable for all $a>0$; the minimal two-dimensional catenoid in $\mathbb R^3$ is $a$-unstable for all $a>0$ and the hyperbolic plane is $a$-stable for $a\le \frac{1}{8}$ and $a$-unstable for $a\ge\frac{1}{8}$, because the hyperbolic space has Gaussian curvature equal to $-1$ and its first eigenvalue of the Laplacian is $-\frac{1}{4}$. Still in dimension two, let $0<a\le\frac{1}{4}$ be a real number, the general idea is that, for a manifold $M$ of negative Gaussian curvature to be $a$-stable, it is necessary that $M$ has a certain volume growth of geodesic balls. This growth would be of increasing order in values of $a$. For example, the plane with the metric $dr^2+r^2(\log (r+e))^2d\theta^2$ has negative Gaussian curvature, it is $a$-unstable for all $a>0$ and the geodesic balls of radio $R$ have area in order of $R^2\log(R)$. The hyperbolic plane is $\frac{1}{4}$-stable and the geodesic balls of radios $R$ have area $2\pi(\cosh^2(R)-1)$. Bérard and Castillon \autocite{berard2014inverse} found the following pattern:

\noindent\textbf{Theorem} [Bérard-Castillon 2010]\hypertarget{1,99}{} \textit {Let $(M,g)$ be a complete noncompact Riemannian surface and
let $W$ be a locally integrable function on $M$, with $W+$ integrable. Assume that the
operator $\Delta + aK + W$ is nonnegative on $M$ and that}

\noindent (i) $a \in (\frac{1}{4},\infty)$, \textit{or}

\noindent (ii) $a =\frac{1}{4}$\textit{, and} $(M,g)$ \textit{has subexponential volume growth, or}

\noindent(iii) $a \in (0, \frac{1}{4})$, and $(M,g)$ has $k_a$-\textit{subpolynomial volume growth, with} $k_a =
\frac{2+4a}{1-4a}$.

\noindent \textit{Then:}

\noindent \textit{(A) The surface $(M,g)$ has finite topology and at most quadratic volume growth.
In particular, $(M,g)$ is conformally equivalent to a closed Riemannian surface
with finitely many points removed.}

\noindent \textit{(B) The function $W$ is integrable on $(M,g)$, and}
$$0 \le 2\pi a \chi(M) + \int_MW.$$

\noindent \textit{(C) If $2\pi a \chi(M)+\int_MW=0$, then $(M,g)$ has subquadratic volume growth,
and $aK + W \equiv 0$ a.e. on the surface $M$.}

This theorem guarantee that if $\alpha\ge 1$ and $(M^2,g)$  has an $(\alpha+1)$ polynomial volume growth of the area, that is, there exists $0<C_1\le C_2<\infty$ such that $C_1\le\frac{|B_R(p)|}{R^{\alpha+1}}\le C_2$ for all $R>0$, then $M$ cannot be $\beta$-stable for all $\beta>\frac{\alpha-1}{8\alpha}$. On the other hand, there are examples of $\frac{\alpha-1}{8\alpha}$-stable manifolds satisfying $|B_R|=O(R^{\alpha+1})$, as exemplified in work \autocite{berard2014inverse} by the same author and in \Cref{cit3.23}. Another important work to be mentioned regarding these operators is that of Espinar and Rosenberg in \autocite{espinar2011colding}.

Suppose that $M$ is a complete bi-dimensional Riemannian manifold with nonpositive Gaussian curvature. To check whether $M$ is $a$-stable, as mentioned in \autocite{kawai1988operator}, we only need to check if the universal cover of $M$ is $a$-stable. Hence, we assume that $M$ is diffeomorphic to $\mathbb R^2$. In this case, the metric $ds^2$ of $M$ can be written as $ds^2= dr^2+\rho(r,\theta)^2d\theta^2$ using polar coordinates centered at any given point.

\section{Case \texorpdfstring{$n\ge 3$}{p} with warped-product metric}

A natural question is whether the results of the previous section hold true for higher dimensions. Answering this question, it is not possible to obtain a general result that generalizes these in dimension two. In this work, we show that some of these results can be generalized to a specific class of Riemannian manifolds, namely to  warped products of an infinite interval
of the line ($[0, \infty)$ or $(-\infty,\infty)=\mathbb R$) by one compact manifolds. Below we define each of these terms.

\begin{definicao} 
Let $(B^m,g_1)$ and $(F^n,g_2)$ be Riemannian manifolds, and $\rho:B\rightarrow\mathbb R$ be a smooth function. We define the \textbf{warped product of $B$ and $F$ with the warping function $\rho$}, denoted by $B\times_{\rho} F$ as the product manifold $B\times F$ with the metric $g$ defined by$$g=g_1+\rho^2g_2.$$
\end{definicao}

We study the stability of the operator $L_a=\Delta-aS$ when $M$ is described as a warped product related to the order of growth of the warping functions. We study for a given growth rate of the warping function, for which the values of $a$ $L_a$ is stable and for which it is unstable. Many results regarding the instability of $L_a$, 
as we will see, reduce to analyze the instability of an end $E$ of $M$, because the growth rate of the warping function determines whether $M$ can or cannot be $a$-stable.

First, we present a result that determines an stability interval of $n$-dimensional warped products.

\begin{teorema}\label{3.1}
    Let $M^n=I\times_{\rho} F^{n-1}$ be a warped product manifold with $n\ge 3$, $I\subset\mathbb R$ and $F$ a compact manifold with nonnegative scalar curvature, then $M$ is $a$-stable for $0\le a\le \frac{n-2}{4(n-1)}$. Consequently, the Yamabe operator $Y\coloneqq\frac{n-1}{4(n-2)}\Delta-S$ is nonnegative  on $M$, being positive definite if $S(F)>0$.
\end{teorema}

The next theorems states that for certain values of $a$ and under certain conditions of the warping function, in some cases, we can guarantee that a (single) warped product is $a$-stable or $a$-unstable. The following theorem resembles the  \hyperlink{Kawai}{Kawai's Theorem} \autocite{kawai1988operator}. 

\begin{teorema}\label{cit3.3}  Let  $M=[0,\infty)\times_{\rho} F$ be a warped product without boundary such that $\rho$ satisfies $\rho''(r)> 0$ and $\lim_{r\rightarrow\infty}\rho'(r)=\infty$. Then the end of $M$ (and therefore, $M$) is $a$-unstable for all $a>\frac{n-1}{4n}$. If $F$ has nonpositive total scalar curvature $S(F)=\int_FS_F$ the hypothesis $\lim_{r\rightarrow\infty}\rho'(r)=\infty$ can be removed.\end{teorema}

\begin{teorema}\label{cit3.4}  Let  $M=I\times_{\rho}F$ be a smooth warped product such that there exists $R_0$ such that  one of the following two situations occurs:

   \noindent (i) $\rho''(r)> 0$ for all $r>R_0$ and $\lim_{r\rightarrow\infty}\rho'(r)=\infty$.

    \noindent (ii) $\rho''(r)> 0$ for all $r<R_0$ and $\lim_{r\rightarrow-\infty}\rho'(r)=-\infty$. 
    
    \noindent Then the corresponding end of $M$ in (i) or (ii) (and therefore, $M$) is $a$-unstable for all $a>\frac{n-1}{4n}$.\end{teorema}

\begin{teorema}\label{cit3.5} Let  $M=I\times_{\rho} \mathbb S^{n-1}$ be a warped product such that $\rho$ satisfies  $|\rho'(r)|\le C$. Then $M$ is $a$-stable for all $0\le a\le \frac{(C^2+1)(n-2)}{4C^2(n-1)}$.\end{teorema}

In the case $n=2$, we saw in \Cref{888} that if in $M$ $C_1\le\frac{|B_R|}{R^{\alpha+1}}\le C_2$, then $M$ is $a$-unstable for all $a>\frac{\alpha-1}{8\alpha}$. The condition $C_1\le \frac{|B_R|}{R^{\alpha+1}}\le C_2$ can be exchanged, in the case of warped products, as $C_1r^{\alpha}\le \rho(r)\le C_2r^{\alpha}$, where $\rho$ is the warping function. As we will see, in a higher dimension and for warped products, this result follows with the hypothesis $C_1r^{\alpha}\le \rho(r)\le C_2r^{\alpha}$, $\alpha>1$, where we will conclude that $M$ is $a$-unstable for all $a>\frac{(n\alpha-\alpha-1)^2}{4\alpha(n-1)(n\alpha-2)}$. Thus obtain an "interval of instability" $\left(\frac{(n\alpha-\alpha-1)^2}{4\alpha(n-1)(n\alpha-2)},\infty\right)$. The case covered in the following theorem associates the polynomial growth of $\rho$ with an "instability interval". 

In the following theorems, we will study the cases where $\rho$ has a growth of polynomial order, focusing on the values of $a$ in the interval $(\frac{n-2}{4(n-1)},\frac{n-1}{4n}]$ to obtain a relationship between $a$-stability and the "degree" of a polynomial which has a growth comparable to the growth of $\rho$. 

\begin{teorema}\label{cit3.8}Suppose $\alpha>1$ and $M=I\times_{\rho}F$ with a metric of the form $g=dr^2+\rho(r)^2g_F$. Suppose that there exist positive constants $C_1$ and $C_2$ such that $C_1r^{\alpha}\le \rho(r)\le C_2r^{\alpha}$ for all $r>1$, then $M$ is $a$-unstable for all $a>\frac{(n\alpha-\alpha-1)^2}{4\alpha(n-1)(n\alpha-2)}.$\end{teorema}

The value $\frac{(n\alpha-\alpha-1)^2)}{4\alpha(n-1)(n\alpha-2)}$ in \Cref{cit3.8} is the best possible because the warped product $M=[0,\infty)\times_{\rho}\mathbb S^{n-1}$ when $\rho(r)=r^{\alpha}$ is $\frac{(n\alpha-\alpha-1)^2}{4\alpha(n-1)(n\alpha-2)}$-stable (see \Cref{cit3.23}). If $\rho$ has $\alpha$-subpolynomial grown, there is a similar version of the \hyperlink{1,99}{Bérard-Castillon Theorem}:

\begin{teorema}\label{cit3.11}Let $M=I\times_{\rho}F$ be a warped product, where $F$ is a compact manifold and $\rho(r)=r^{\alpha}\xi(r)$, where $\alpha\ge 1$ and $\xi$ satisfies $\xi(r)\rightarrow 0$ when $r\rightarrow\infty$. Suppose that $\xi$ is a non-increasing function and satisfies $r\xi(r)^{-1}\xi'(r)\ge -(\alpha-1)$ for all $r$ large. Also suppose that $M$ is $a$-stable for $a= \frac{(n\alpha-\alpha-1)^2}{4\alpha(n-1)(n\alpha-2)}$. Then $\rho$ has linear growth, that is, $\rho(r)\le Cr$ for some $C\ge 0$. In particular, $M$ has a polynomial volume growth on order of $R^n$, that is, for each $p\in M$, $Vol(B_R(p))\le C_1R^n$ for some $C_1>0$.\end{teorema}

Consider the application \begin{equation*}\begin{split}h:\left(1,\infty\right) & \longrightarrow \left(\frac{n-2}{4(n-1)},\frac{n-1}{4n}\right)\\ \zeta & \longmapsto \frac{(n\zeta-\zeta-1)^2}{4\zeta(n-1)(n\zeta-2)}.\\ \end{split}\end{equation*}

\noindent Using elementary analysis tools, it is possible to show that $h$ is an increasing function, $$\lim_{\zeta\rightarrow 1}h(\zeta)=\frac{n-2}{4(n-1)}\text{\hspace{4mm} and\hspace{4mm}}\lim_{\zeta\rightarrow \infty}h(\zeta)=\frac{n-1}{4n}.$$ 
Therefore, there is a continuous and monotonous relationship between the degree of polynomial growth of $\rho$ and the maximum possible value of $a$ for which $M$ can be $a$-stable.

The next theorem is an improved version of the \Cref{cit3.3} for when we assume that $\rho$ has a polynomial growth.

\begin{teorema}\label{xiiso(1)}Let $M=I\times_{\rho}F$ be a warped product, where $F$ is a compact manifold of dimension at least two and $\rho(r)=r^{\alpha}\xi(r)$, where $\alpha$ and $\xi$ satisfy $$\alpha=\inf\left\{\gamma;\hspace{2mm}\lim_{r\rightarrow\infty}\rho(r)r^{-\gamma}=0\right\}.$$ Suppose that $M$ is $a$-stable for some $a>\frac{n-1}{4n}$. Then there exist a positive constant $C$ such that $$\liminf_{R\rightarrow\infty}(\log R)^{-1}\int_1^Rr\rho(r)^{-2}dr\ge C.$$ Furthermore, if $S(F)\le 0$, then $\alpha<\frac{n}{2n-2}$. In particular, $M$ has $(1+\frac{n}{2})$-subpolynomial volume growth, that is, $$\lim_{r\rightarrow\infty}r^{-(1+\frac{n}{2})}Vol(B_r(p))=0$$ for all $p\in M$.\end{teorema}

The following theorem shows what can happen if the function $\xi$ as in the previous theorem has large variation.

\begin{teorema}\label{th3.1.8}Let $M=I\times_{\rho}F$, where $\rho(r)=r^{\alpha}\xi(r)$, where $\alpha\ge 1$ and $\xi$ satisfy $$\alpha=\inf\left\{\gamma;\hspace{2mm}\lim_{r\rightarrow\infty}\rho(r)r^{-\gamma}=0\right\}$$ with the additional condition $C_1\le\xi(r)\le C_2$ for all $r\ge 0$ for some constants $C_1,C_2>0$. Suppose that $$\lim_{T\rightarrow\infty}\frac{log T}{\int_{R_0}^Tr\xi(r)^{-2}\xi'(r)^2dr}=0,$$ where $R_0$ is a fixed positive number. Then $M$ is $a$-unstable for all $a>\frac{n-2}{4(n-1)}$.\end{teorema}

In \Cref{cit3.8} we assumed that $\alpha>1$. The following results exemplify what can happen if $\alpha=1$, serving as an object of comparison with the previous result in which the assumption $\alpha>1$ was required.

\begin{teorema}\label{cit3.12}  Let $M=I\times_{\rho}\mathbb S^{n-1}$ be a warped product, where $I=\mathbb R$ or $I=[0,\infty)$. Suppose that $C_1r\le \rho(r)\le C_2r$ $\forall r\ge 1$ for the positive constants $C_1,C_2$. If $C_1>1$, then (the ends of) $M$ is $a$-unstable for all $a>\frac{n-2}{4(n-1)(1-C_1^{-2})}$.\end{teorema}

\begin{teorema} \label{cit3.13} Let $M=I\times_{\rho}\mathbb S^{n-1}$ be a warped product, where $I=\mathbb R$ or $I=[0,\infty])$. Suppose that $C_1r\le \rho(r)\le C_2r$ $\forall r\ge 1$ for the positive constants $C_1,C_2$. If $C_2<1$, then (the ends of) $M$ is $a$-unstable for all $a<\frac{n-2}{4(n-1)(1-C_2^{-2})}$.\end{teorema}

Thus, based on these results, we obtain the following diagram, relating the $a$-stability of a smooth manifold $M=I\times_{\rho}\mathbb S^{n-1}$ with the possible growth characteristics of the $\rho$ function
:

\

\begin{center}{\includegraphics[scale=0.38306]{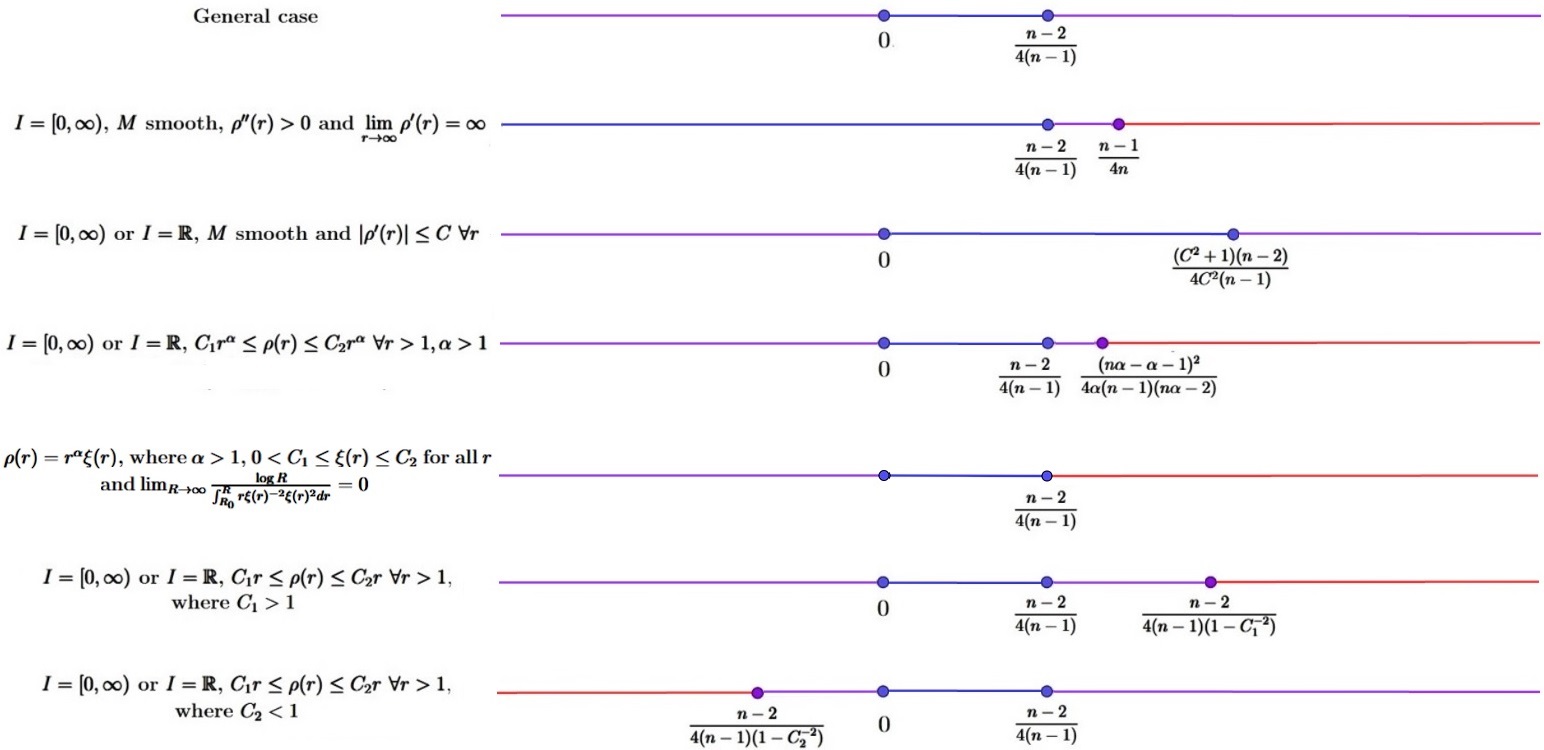}}\end{center}

\crule[blue!60]{3mm}{3mm} \scriptsize{$a$-stable}

\crule[red]{3mm}{3mm} $a$-unstable

\crule[red!60!blue!80!white!80!]{3mm}{3mm} can be $a$-stable or $a$-unstable

\normalsize

\section{Operator \texorpdfstring{$\Delta-aS$}{d} on Warped Products}\label{C3}

When we make a warped product of two smooth manifolds with a positive warping function $\rho$, the resulting manifold is necessarily a smooth manifold. If we require $\rho$ to be nonnegative, the points in the base such that $\rho$ vanishes represent a degenerate fiber, in the sense that becomes one point. The resulting manifold of a warped product of a closed interval of $\mathbb R$ with a manifold $F$ can result in a smooth manifold. 

Consider an elliptic operator of the form $L=\Delta-q$ on a warped product $M=I\times_{\rho}F$, where $\Delta$ is the Laplace-Beltrami operator of $M$, $q:M\rightarrow \mathbb R$ is a smooth function, $I$ is a closed interval of $\mathbb R$ (which we can essentially reduce to three cases: $\mathbb R$, $[0,\infty)$ and $[0,1]$), $F$ is a compact manifold and $\rho$ is a nonnegative warping function. Unless otherwise mentioned, we will always assume that $M$ is without boundary, which for us will be the same as requiring that $\rho$ vanishes in $\partial I$. The operator $L$ is said to be \textbf{stable} if $$\int_M-fLf\ge 0$$ for all $f\in C_c^{\infty}(U)$, where $U$ is the set of regular points of $M$, that is, the largest subset of $M$ which is a manifold. In particular, if $V=\{r\in I;\rho(r)>0\}$, then $V\times_{\rho}F\subset U$.
The previous inequality is equivalent to $$\int_M||\nabla f||^2+qf^2\ge 0$$ for all  $f\in C_c^{\infty}(U)$. 

Let $F$ be an $(n-1)$-dimensional compact manifold with $n\ge 3$ and $M=I\times_{\rho}F$, such that the metric of $M$ is expressed in the form $$g=dr^2+\rho(r)^2g_{F},$$

\noindent where $\rho$ is a positive real function in the interior of $I$. We have $M$ compact if and only if $I$ is a finite closed interval. $M$ there will be no boundary if $\rho$ is zero at the endpoints of $I$. $M$ will be complete and noncompact if and only if $I$ is an unbounded interval. Up to a linear variable change in $I$, there are two types of closed and unbounded intervals in $\mathbb R$: $[0,\infty)$ and $\mathbb R$ itself.

Consider the case in which $I=[0,\infty)$ and $\rho> 0$ in $(0,\infty)$. Under this condition, $M$ has one end, $M$ is a manifold with boundary if $\rho(0)>0$ and $M$ is smooth without boundary if and only if $\rho(0)=0$, $F$ is a round sphere of radius $R$ and $\rho'(0)=\frac{1}{R}$, $R>0$. In the case $I=\mathbb R$, $M$ will have two ends and will be smooth if and only if $\rho$ is positive everywhere. In all the cases, the volume element of $M$ is $\rho(r)^{n-1}drdA$, where $dA$ is the area element of $\mathbb S^{n-1}$.

We studied the stability of the operator $L=\Delta-aS$ on $M$ and sometimes at an end of $M$. We say that $L$ is \textit{stable on an end $E$ of $M$}  if the operator $L$ restricted to a representative of $E$ of the form $((R,\infty)\times_{\rho} F,g)$ (or $((-\infty,-R)\times F,g)$) for all $R$ sufficiently large is stable. If $E$ is not stable, we say that $E$ is \textit{unstable}. It is obvious that an end $E$ be unstable implies that $M$ is unstable. 

For warped products, the scalar curvature can be written as \begin{equation*}\label{scalarwarpedformula}S_M=\frac{S_F}{\rho(r)^2}-(n-1)\left(2\frac{\rho''(r)}{\rho(r)}+(n-2)\frac{\rho'(r)^2}{\rho(r)^2}\right).\end{equation*}

For $f$ of compact support in $M$, only dependent on $r$ and support contained in the interval $[b,c]$, we have:

\noindent (i) if $I=\mathbb R$ ,then $f(b)=f(c)=0$; 

\noindent (ii) if $I=[0,\infty)$ and $b\ne 0$, then $f(b)=0$. If $b=0$, we only guarantee $f(c)=0$, but if  we demand $M$ without boundary, then $\rho(b)=\rho(0)=0$. 

\noindent Hence, in the case $M$ without boundary and $f$ of compact support, we have $$\rho(b)^{n-2}f(b)=\rho(c)^{n-2}f(c)=0.$$ Therefore, we have: 
\begin{equation}\label{eq1}\begin{split}\int_M-fL_af = & \int_M|\nabla f|^2+aSf^2
\\ = & \int_b^c\int_F\left\{|\nabla f|^2+a\left[-(n-1)\left(2\frac{\rho''(r)}{\rho(r)}+(n-2)\frac{\rho'(r)^2}{\rho(r)^2}\right)+\frac{S_F}{\rho(r)^2}\right]f^2\right\}dAdr. \end{split}\end{equation}

First, let us look at a result that will make our calculations of whether $M$ is or is not stable. We will use principally to prove stability in a warped product of an interval with round spheres and an interval with compact manifolds of null scalar curvature.

\begin{proposicao}\label{cit3.1}
    Let $M^n=I\times_{\rho}F^{n-1}$ be a warped product, with $n\ge 3$. Suppose that $F$ has constant scalar curvature  and $$-\int_MfL_af\ge 0$$ 
    
    \noindent for all $f$ dependent only on $I$, that is, for all $f$ constants on the fibers. Then $M$ is $a$-stable.
\end{proposicao}

\noindent \textbf{Proof}. \noindent If $f$ is a function of compact support in $M$, then $\pi_1(supp (f))$ has compact support in $I$. If $f$ depends only on $I$, by (\ref{eq1}), we have
\begin{equation}\label{2.2} \int_I \left\{f_r^2+a\left[ -(n-1)\left(2\frac{\rho''(r)}{\rho(r)}+(n-2)\frac{\rho'(r)^2}{\rho(r)^2}\right)+\frac{S(F)}{A(F)\rho(r)^2}\right]f^2\right\}dr\ge 0,\end{equation}

\noindent where $S(F)$ is the total scalar curvature of $F$ and $A(F)$ is the area of $F$. Let $f$ be any function of compact support on $M$, using Fubini's theorem and the fact that $|\nabla f|^2\ge f_r^2$, can be obtained by (\ref{eq1}):

\vspace{10pt}

\noindent$\displaystyle\int_M-fL_af$

\noindent$\displaystyle=\int_F\int_b^c\left\{|\nabla f|^2+a\left[-(n-1)\left(2\frac{\rho''(r)}{\rho(r)}+(n-2)\frac{\rho'(r)^2}{\rho(r)^2}\right)+\frac{S_F(q)}{\rho(r)^2}\right]f(r,q)^2\right\}drdA$

\noindent$\displaystyle\ge \int_F\int_b^c\left\{f_r(r,q)^2+a\left[-(n-1)\left(2\frac{\rho''(r)}{\rho(r)}+(n-2)\frac{\rho'(r)^2}{\rho(r)^2}\right)+\frac{S_F(q)}{\rho(r)^2}\right]f(r,q)^2\right\}drdA$,

\vspace{10pt}

\noindent where $[b,c]=\pi_1(supp(f))$. Since $f_r(r,q)$ and $f(r,q)$  depend only of $f$ on $r$ in each leaf $I\times \{q\}$ and $\frac{S_F}{\rho(r)^2}$ is constant on $F$, the last expression is equal to
$$\int_F\int_b^c\left\{f_r(r,q)^2+a\left[-(n-1)\left(2\frac{\rho''(r)}{\rho(r)}+(n-2)\frac{\rho'(r)^2}{\rho(r)^2}\right)+\frac{S(F)}{A(F)\rho(r)^2}\right]f(r,q)^2\right\}drdA,$$

\noindent that is nonnegative because (\ref{2.2}). \\ \noindent$\blacksquare$

The \Cref{cit3.1} shows that, with the hypothesis that $F$ has constant scalar curvature, we simply take real functions $f$ as a test function (that is, dependent only on the first coordinate) to demonstrate whether $M$ is stable.

Continuing the calculation in (\ref{eq1}), we have, for $f$ of compact support dependent only on $r$:
\begin{equation*}\begin{split}
\int_M-fL_af= & A(F)\int_b^cf_r^2\rho(r)^{n-1}dr-2a(n-1)\int_b^c\rho''(r)\rho(r)^{n-2}dr \\ &+(n-2)\int_b^c\rho'(r)^2\rho(r)^{n-3})f^2dr +aS(F)\int_b^c\rho(r)^{n-3}f^2dr \\  =& A(F)\int_b^cf_r^2\rho(r)^{n-1}dr+a(n-1)(n-2)A(F)\int_b^c\rho'(r)^2\rho(r)^{n-3}f^2dr\\ & +4a(n-1)A(F)\int_b^c\rho'(r)\rho(r)^{n-2}ff_rdr+aS(F)\int_b^c\rho(r)^{n-3}f^2dr,\end{split}\end{equation*}

\noindent where $A(F)$ is the area of $F$ and $S(F)$ is the total scalar curvature of $F$. Therefore
\begin{equation}\label{2.3}\begin{split}\frac{1}{A(F)}\int_M-fL_af = & \int_b^cf_r^2\rho(r)^{n-1}dr+a(n-1)(n-2)\int_b^c\rho'(r)^2\rho(r)^{n-3}f^2dr+ \\
& +4a(n-1)\int_b^c\rho'(r)\rho(r)^{n-2}ff_rdr+a\frac{S(F)}{A(F)}\int_b^c\rho(r)^{n-3}f^2dr.\end{split}\end{equation}

\noindent From equation (\ref{2.3}), we can find an interval on the line for which the operator $L_a$ is always stable in any warped product of form $I\times_{\rho}F$, whenever $F$ has nonnegative total scalar curvature.

Assuming that $\rho$ has a polynomial growth greater than the Euclidean, that is, there exists $1<\zeta<\infty$ such that $\lim_{r\rightarrow\infty}\rho(r)r^{-\zeta}=0$ and $\lim_{r\rightarrow\infty}\rho(r)r^{-1}=+\infty$. Here, assuming $M=I\times_{\rho}F$ with $I=[0,\infty)$ or $I=\mathbb R$, we are analyzing the end of $M$ corresponding to the positive end of the real line $\mathbb R$. The analysis of the negative end of $\mathbb R$, when $I=\mathbb R$, is analogous, taking the application $\tilde{\rho}$ defined by $\tilde{\rho}(r)=\rho(-r)$ and analyzing $\tilde{\rho}$ on the positive end of $\mathbb R$.  Under these conditions, there exists a smooth nonnegative real function $\xi:M\rightarrow\mathbb R$ and a real number $\alpha\ge b1$ such that for $r\ge 1$, \begin{equation}\label{3.6}\rho(r)=r^{\alpha}\xi(r)\text{, where }\alpha=\inf\left\{\gamma;\hspace{2mm}\lim_{r\rightarrow\infty}\rho(r)r^{-\gamma}=0\right\}.\end{equation} Many situations in relation to $\xi$ can occur, among these, we can have $\lim_{r\rightarrow\infty}\xi(r)=\infty$ (exemplified by $\xi(r)=\log(r^2+2)^k$, $k>0$), we can have $\lim_{r\rightarrow\infty}\xi(r)=0$ (exemplified by $\xi(r)=\log(r^2+2)^k$, $k<0$) and the case where there exist positive constants $C_1$ and $C_2$ such that $C_1\le\xi(r)\le C_2$ as when $\xi$ is constant. If $\lim_{r\rightarrow\infty}\xi(r)=\infty$, we have $\lim_{r\rightarrow\infty}\xi(r)r^{-\gamma}=0$ for all $\gamma>0$ and if $\lim_{r\rightarrow\infty}\xi(r)=0$, we have $\lim_{r\rightarrow\infty}\xi(r)r^{\gamma}=\infty$ for all $\gamma>0$. Being $f$ of compact support only dependent on $r$, assume that the support of $f$ is of form $[b,c]\times F$, with $b,c\neq 0$. 
Being $\rho(r)=r^{\alpha}\xi(r)$, we have, for $r\neq 0$, $\rho'(r)=\alpha r^{\alpha-1}\xi(r)+r^{\alpha}\xi'(r)$ and $$\rho'(r)^2=\alpha^2r^{2\alpha-2}\xi(r)^2+2\alpha r^{2\alpha-1}\xi(r)\xi'(r)+r^{2\alpha}\xi'(r)^2.$$ Then (\ref{2.3}) becomes
\begin{equation}\label{2.4}\begin{split}\frac{1}{A(F)}\int_M-fL_af =& \int_b^cf_r^2r^{n\alpha -\alpha}\xi(r)^{n-1}dr+a\alpha^2(n-1)(n-2)\int_b^cr^{n\alpha-\alpha-2}\xi(r)^{n-1}f^2dr\\ &+2a\alpha(n-1)(n-2)\int_b^cr^{n\alpha-\alpha-1}\xi(r)^{n-2}\xi'(r)f^2dr \\ &+a(n-1)(n-2)\int_b^cr^{n\alpha-\alpha}\xi(r)^{n-3}\xi'(r)^2f^2dr \\ &+4a\alpha(n-1)\int_b^cr^{n\alpha-\alpha-1}\xi(r)^{n-1}ff_rdr \\ &+4a(n-1)\int_b^cr^{n\alpha-\alpha}\xi(r)^{n-2}\xi'(r)ff_rdr \\ &+a\frac{S(F)}{A(F)}\int_b^cr^{n\alpha-3\alpha}\xi(r)^{n-3}f^2dr.\end{split}\end{equation}

\section{Proofs of the Theorems}\label{S3.3}

\noindent\textbf{Proof of \Cref{3.1}}. By (\ref{eq1}) taking $f:M\rightarrow\mathbb R$ of compact support and $\pi_1(Supp f)\subset [b,c]$, using Fubini's theorem, integration by parts, that $S_F\ge 0$ and the fact $$\rho(b)^{n-2}f(b)=\rho(c)^{n-2}f(c)=0,$$ we obtain:
\begin{equation}\label{eq6}\begin{split}\int_M-fL_af= & \int_M|\nabla f|^2+\int_F\left[a(n-1)(n-2)\int_b^c\rho'(r)^2\rho(r)^{n-3}f(r,q)^2dr\right. \\& \left.+4a(n-1)\int_b^c\rho(r)^{n-2}\rho'(r)f(r,q)f_r(r,q)dr+\int_b^c\frac{S_F(q)}{\rho(r)^2}dr\right]dA \\ \ge & \int_F\left[\int_b^cf_r(r,q)^2\rho(r)^{n-1}dr+a(n-1)(n-2)\int_b^c\rho(r)^{n-3}\rho'(r)^2f(r,q)^2dr \right.\\  & \left.+4a(n-1)\int_b^c\rho(r)^{n-2}\rho'(r)f(r,q)f_r(r,q)dr\right]dA.\end{split}\end{equation}

\noindent Using the AM-GM inequality, we have, for each $q\in F$:

\noindent$\displaystyle-4a(n-1)\int_b^c\rho(r)^{n-2}\rho'(r)f(r,q)f_r(r,q)dr$

\noindent$\displaystyle=-4a(n-1)\int_b^cf_r(r,q)\rho(r)^{\frac{n-1}{2}}\rho(r)^{\frac{n-3}{2}}\rho'(r)f(r,q)dr$

\noindent$\displaystyle \le \int_b^cf_r(r,q)^2\rho(r)^{n-1}dr  +4a^2(n-1)^2\int_b^c\rho(r)^{n-3}\rho'(r)^2f(r,q)^2dr$.

\noindent Then
\begin{equation}\label{2.04}\begin{split} &\int_b^cf_r(r,q)^2\rho(r)^{n-1}dr+4a^2(n-1)^2\int_b^c\rho(r)^{n-3}\rho'(r)^2f(r,q)^2dr \\ &+4a(n-1)\int_b^c\rho(r)^{n-2}\rho'(r)f(r,q)f_r(r,q)dr\ge 0.\end{split}\end{equation}

\noindent So to conclude the proof, we just need to show that
\begin{equation*}\begin{split} & a(n-1)(n-2)\ge 4a^2(n-1)^2 \\ \Longleftrightarrow & a\ge 0 \text{\hspace{3mm}and\hspace{3mm}} n-2\ge 4a(n-1)\\ \Longleftrightarrow & a\ge 0 \text{\hspace{3mm}and\hspace{3mm}} a\le\dfrac{n-2}{4(n-1)}.\end{split}\end{equation*}

\noindent This ensures the $a$-stability of $M$ for all $a\in [0,\frac{n-2}{4(n-1)}]$. Taking $a=\frac{n-2}{4(n-1)}$, then the Yamabe operator $Y$ is nonnegative. 
If $S(F)>0$, then the existence of the term $\int_b^c\frac{S_F(q)}{\rho(r)^2}dr$ in the calculation (\ref{eq6}) implies the positively of $Y$. If $I$ is unbounded and $S(F)=0$, which in our hypotheses means that the scalar curvature of $F$ is identically zero, suppose that $f\in C_c^{\infty}(M)$ satisfies $$-\int_M fY(f)=0,$$ 
\noindent then the equality in (\ref{2.04}) occurs, it implies $$\int_b^c\left[f_r(r,q)\rho(r)^{\frac{n-1}{2}}+\frac{n-2}{2}\rho(r)^{\frac{n-3}{2}}\rho'(r)f(r,q)\right]^2dr=0$$

\noindent for all $q\in F$, this implies \begin{equation*}\begin{split}f_r(r,q)=-\frac{n-2}{2}\rho
(r)^{-1}\rho'(r)f(r,q)\Longrightarrow & \log(f(r,q))=-\frac{n-2}{2}\log(\rho(r))+C(q) \\ \Longrightarrow & f(r,q)=C(q)\rho(r)^{-\frac{n-2}{2}},\end{split}\end{equation*} that has compact support if and only if $C(q)=0$ for all $q\in F$, that is, if and only if $f\equiv 0$. 

\noindent$\blacksquare$

\noindent\textbf{Proof of \Cref{cit3.3}}. First, since $M$ is without boundary, $\rho(0)=0$ and by hypothesis $\rho''(r)>0$ for all $r$, we have $r\rho'(r)>\rho(r)$ for $r>0$, because for $r=0$ the equality holds and the application $r\mapsto r\rho'(r)-\rho(r)$ has as derivative the application $r\mapsto r\rho''(r)>0$. Therefore, $r\rho(r)^{-1}\rho'(r)>1$ for all $r>0$. Furthermore, the application $r\mapsto r\rho(r)^{-1}$ is decreasing because its derivative is $$r\mapsto \rho(r)^{-1}-r\rho(r)^{-2}\rho'(r)=\rho(r)^{-1}(1-r\rho(r)^{-1}\rho'(r))<0.$$ Rewriting (\ref{2.3}) we have, for $f=f(r)$ of compact support contained in $[b,c]$: 
\begin{equation}\label{9}\begin{split}\frac{1}{A(F)}\int_M-fL_af= & \int_b^cf_r^2\rho(r)^{n-1}dr+a(n-1)(n-2)\int_b^c\rho'(r)^2\rho(r)^{n-3}f^2dr \\ &
+4a(n-1)\int_b^c\rho'(r)\rho(r)^{n-2}ff_rdr+a\frac{S(F)}{A(F)}\int_b^c\rho(r)^{n-3}f^2dr.\end{split}\end{equation}

\noindent Take $$f(r)=f_{Q,R}(r)=\begin{cases} \hspace{31.5mm}0\hspace{31.5mm}\text{if}\hspace{3mm}0\ge \frac{Q}{2} \\ 2\rho(r)^{-\frac{1}{2}(n-1)}Q^{-\frac{1}{2}}r-\rho(r)^{-\frac{1}{2}(n-1)}Q^{\frac{1}{2}}\hspace{2mm}\text{if}\hspace{3mm}\frac{Q}{2}\le r\le Q; \\ \hspace{19.5mm}\rho(r)^{-\frac{1}{2}(n-1)}r^{\frac{1}{2}},\hspace{19.5mm}\text{if}\hspace{3mm}Q\le r\le R; \\ C(\rho(r)^{-\frac{1}{2}(n-1)}r^{\frac{1}{2}}-\rho(T)^{-\frac{1}{2}(n-1)}T^{\frac{1}{2}})  \hspace{3mm}\text{if}\hspace{3mm}R\le r\le T; \\ \hspace{31.5mm}0\hspace{31.5mm}\text{if}\hspace{3mm}r\ge T, \end{cases}$$

\noindent where $T$ is such that $\rho(T)=2\rho(R)$ and $C$ is such that becomes $f$ a continuous function. Assume $Q$ is large such that $r\rho(r)^{-1}\rho'(r)>1$ for all $r>Q$. We fix $Q$, then the integration of the expression (\ref{9}) in the interval $[\frac{Q}{2},Q]$ is fixed. In the interval $[Q,R]$ it becomes, using that $f_r=-\frac{1}{2}(n-1)\rho(r)^{-\frac{1}{2}(n-1)-1}\rho'(r)r^{\frac{1}{2}}+\frac{1}{2}\rho(r)^{-\frac{1}{2}(n-1)}r^{-\frac{1}{2}}$:

\vspace{10pt}

\noindent \hspace{4mm}$\displaystyle\frac{1}{4}(n-1)^2\int_Q^Rr\rho(r)^{-2}\rho'(r)^2dr-\frac{1}{2}(n-1)\int_Q^R\rho(r)^{-1}\rho'(r)dr+\frac{1}{4}\int_Q^Rr^{-1}dr$

\noindent \hspace{4mm}$\displaystyle+a(n-1)(n-2)\int_Q^Rr\rho(r)^{-2}\rho'(r)^2dr+4a(n-1)\int_Q^R-\frac{1}{2}(n-1)r\rho(r)^{-2}\rho'(r)^2dr$

\noindent\hspace{4mm}$\displaystyle+4a(n-1)\int_Q^R\frac{1}{2}\rho(r)^{-1}\rho'(r)dr+a\frac{S(F)}{A(F)}\int_Q^Rr\rho(r)^{-2}dr$

\noindent$\displaystyle=\frac{1}{4}(n-1)(n-1-4an)\int_Q^Rr\rho(r)^{-2}\rho'(r)^2dr+\left(2a-\frac{1}{2}\right)(n-1)\int_Q^R\rho(r)^{-1}\rho'(r)dr$

\noindent\hspace{4mm}$\displaystyle+\frac{1}{4}\int_Q^Rr^{-1}dr+a\frac{S(F)}{A(F)}\int_Q^Rr\rho(r)^{-2}dr$.
\begin{equation}\label{eq9}\end{equation}
Since $\rho'(r)\rightarrow \infty$ when $r\rightarrow\infty$, the last term of (\hyperlink{9}{9}) is arbitrary small comparing to $$\int_Q^Rr\rho(r)^{-2}\rho'(r)^2dr$$ if we choose $Q$ large. Since $r\rho(r)^{-1}\rho'(r)>1$, we have $\rho(r)^{-1}\rho'(r)>\frac{1}{r}$, implying $r\rho(r)^{-2}\rho'(r)^2>\rho(r)^{-1}\rho'(r)>\frac{1}{r}$. We wish to show that (\hyperlink{9}{9}) is negative and upper limited by a negative constant $C_1$ times $\int_Q^Rr\rho(r)^{-2}\rho'(r)^2dr$. We have $(2a-\frac{1}{2})(n-1)+\frac{1}{4}<0$ if and only if $a<\frac{2n-3}{8(n-1)}$, then this statement is true for $\frac{n-1}{4n}<a\le \frac{2n-3}{8(n-1)}$ taking $C_1=\frac{1}{4}(n-1)(n-1-4an)$. For $a>\frac{2n-3}{8(n-1)}$, we only have to show that $$\frac{1}{4}(n-1)(n-1-4an)+\left(2a-\frac{1}{2}\right)(n-1)+\frac{1}{4}<0,$$

\noindent that is equivalent to
\begin{equation}\label{10}(n-1)\left(-(n-2)a+\frac{n}{4}-\frac{3}{4}\right)<-\frac{1}{4}.\end{equation}

\noindent Since $n\ge 3$, the left-hand side of the previous inequality decreases as a function of $a$. At $a=\frac{2n-3}{8(n-1)}$ the left-hand side of (\ref{10}) becomes $$(n-1)\left(\frac{n}{4}-\frac{(2n-3)(n-2)}{8(n-1)}-\frac{3}{4}\right)=\frac{2n^2-2n-2n^2+7n-6-6n+6}{8}=-\frac{n}{8}<-\frac{1}{4}.$$

\noindent Hence, it follows from the decreasing property that the inequality (\ref{10}) holds for all $a>\frac{2n-3}{8(n-1)}$.  Note that $\int_Q^Rr\rho(r)^{-2}\rho'(r)^2dr\rightarrow\infty$ when $R\rightarrow\infty$, because $r\rho(r)^{-2}\rho'(r)^2>\frac{1}{r}$ for $r>Q$. 
Let us now analyze the expression of $f$ in (\ref{9}) in the interval $[R,T]$. We have $$f_r=-\frac{1}{2}(n-1)\rho(r)^{-\frac{1}{2}(n+1)}\rho'(r)r^{\frac{1}{2}}+\frac{1}{2}\rho(r)^{-\frac{1}{2}(n-1)}r^{-\frac{1}{2}}.$$ Then the expression of $f$ in (\ref{9}) in interval $[R,T]$ is $C^2$ times the expression
\begin{equation}\label{eq11}\begin{split}& \frac{1}{4}(n-1)^2\int_R^Tr\rho(r)^{-2}\rho'(r)^2dr-\frac{1}{2}(n-1)\int_R^T\rho(r)^{-1}\rho'(r)dr+\frac{1}{4}\int_R^Tr^{-1}dr \\ & +a(n-1)(n-2)\int_R^Tr\rho(r)^{-2}\rho'(r)^2dr \\ &-2a(n-1)(n-2)\rho(T)^{-\frac{1}{2}(n-1)}T^{\frac{1}{2}}\int_R^T\rho(r)^{\frac{n-5}{2}}\rho'(r)^2r^{\frac{1}{2}}dr\\ &+a(n-1)(n-2)\rho(T)^{-n+1}T\int_R^T\rho(r)^{n-3}\rho'(r)^2dr-2a(n-1)^2\int_R^Tr\rho(r)^{-2}\rho'(r)^2dr\\ &+2a(n-1)^2\rho(T)^{-\frac{1}{2}(n-1)}T^{\frac{1}{2}}\int_R^Tr^{\frac{1}{2}}\rho(r)^{\frac{n-5}{2}}\rho'(r)^2dr+2a(n-1)\int_R^T\rho(r)^{-1}\rho'(r)dr\\ &-2a(n-1)\rho(T)^{-\frac{1}{2}(n-1)}T^{\frac{1}{2}}\int_R^Tr^{-\frac{1}{2}}\rho(r)^{\frac{n-3}{2}}\rho'(r)dr+a\frac{S(F)}{A(F)}\int_R^Tr\rho(r)^{-2}dr\end{split}\end{equation}

\noindent Being $\rho(T)=2\rho(R)$, by the mean value theorem and using that $\rho''(r)>0$ and $$r\rho(r)^{-1}\rho'(r)>1\Longleftrightarrow r\rho'(r)>\rho(r)$$ for all $r>Q$, there exists $c\in (R,2R)$ such that  $$\rho(2R)-\rho(R)=R\rho'(c)>R\rho'(R)>\rho(R)\Longrightarrow \rho(2R)>2\rho(R).$$ Hence $R<T<2R$. The relation $\rho(T)=2\rho(R)$ with $R<T<2R$ implies that all the terms of (\ref{eq11}) are upper limited by a multiple of $$\int_R^Tr\rho(r)^{-2}\rho'(r)^2dr,$$ because all the expressions are limited by a constant times $\int_R^Tr\rho(r)^{-2}\rho'(r)^2dr$. We observe the expression $\int_R^Tr\rho(r)^{-2}\rho'(r)^2dr$, as a consequence of integration by parts and using that $\rho(T)=2\rho(R)$, $r\rho(r)^{-1}$ is decreasing and $\rho''(r)\ge 0$ for $r>Q$, we have
\begin{equation}\label{eq12}\begin{split}\int_R^Tr\rho(r)^{-2}\rho'(r)^2dr= & \int_R^T\rho(r)^{-1}\rho'(r)dr+\int_R^Tr\rho(r)^{-1}\rho''(r)dr-r\rho(r)^{-1}\rho'(r))|_R^T\\ \le & \log (\rho(T))-\log(\rho(R))+R\rho(R)^{-1}(\rho'(T)-\rho'(R))-T\rho(T)^{-1}\rho'(T)\\ & +R\rho(R)^{-1}\rho'(R)\\ = & \log 2+2R\rho(T)^{-1}\rho'(T)-T\rho(T)^{-1}\rho'(T)\\ \le & \log 2+T\rho(T)^{-1}\rho'(T)\end{split}\end{equation}




\noindent Now, our objective is to show that the expression of the stability operator in $[Q,R]$ can be arbitrarily larger than the expression of the stability operator in $[R,T]$ for $R$ and $T$ appropriated. Let $g(r)=\rho(r)^{-1}\rho'(r)$. We assert that for all $C>0$, there exists $V$ sufficiently large such that \begin{equation}\label{11}\int_Q^Vrg(r)^2dr>CVg(V).\end{equation}

\noindent For this, suppose that $\int_Q^Vrg(r)^2dr\le CVg(V)$ for some constant $C$ and for all $V>Q$, let $h(V)=\int_Q^Vrg(r)^2dr$, then $$h(V)\le C(Vh'(V))^{\frac{1}{2}}\Longrightarrow h(V)^2\le C^2Vh'(V)\Longrightarrow \frac{h'(V)}{h(V)^2}\ge \frac{1}{C^2V}.$$

\noindent Integrating the last inequality from $Q$ to $V$, we obtain $$h(Q)^{-1}-h(V)^{-1}\ge \frac{1}{C^2}(\log V-\log Q).$$ 

\noindent This is a contradiction because we saw that $h(V)=\int_Q^Vr\rho(r)^{-2}\rho'(r)^2dr$ tends to infinity when $V$ tends to infinity, the left-hand side of the last inequality is upper limited by $h(Q)^{-1}$ and the right-hand side tends to infinity when $V$ tends to infinity. 

That said, combining (\ref{eq12}) and (\ref{11}), there is an increasing and unbounded sequence $$\{T_1,T_2,\dots\}$$ such that $$\lim_{j\rightarrow\infty}\frac{\int_Q^{T_j}r\rho(r)^{-2}\rho'(r)^2dr}{\int_{R_j}^{T_j}r\rho(r)^{-2}\rho'(r)^2dr}=\infty,$$

\noindent where $R_j$ is such that $\rho(T_j)=2\rho(R_j)$. It implies $$\lim_{j\rightarrow\infty}\frac{\int_Q^{R_j}r\rho(r)^{-2}\rho'(r)^2dr}{\int_{R_j}^{T_j}r\rho(r)^{-2}\rho'(r)^2dr}=\infty.$$ 

Therefore, the expression of the stability operator in $[Q,R_j]$ is arbitrarily larger than the expression of the stability operator in $[R_j,T_j]$ when $J$ increases to infinity. It finishes the proof.

Note that if $F$ has nonpositive total scalar curvature, the last term of (\hyperlink{9}{9}) and (\ref{eq11}) will be negative and can be disregarded in the respective calculus of instability and we do not need of $\rho'(r)\rightarrow\infty$ when $r\rightarrow\infty$.

\noindent $\blacksquare$

\noindent \textbf{Proof of \Cref{cit3.4}}: Suppose that (i) occurs. The proof is similar to \Cref{cit3.3}, with some modifications. 
We affirm that it continues to be valid that $r\rho'(r)>\rho(r)$ for $r$ sufficiently large (larger than $R_0$). To verify this, because $(r\rho'(r)-\rho(r))'=r\rho''(r)>0$, it follows that $r\rho'(r)-\rho(r)$ increases and we only have to prove that the inequality $r\rho'(r)>\rho(r)$ holds for one value of $r$. Suppose that $r\rho'(r)\le\rho(r)$ for all $r>R_0$, then for $r>R_0$:
\begin{equation*}\begin{split}\frac{\rho'(r)}{\rho(r)}\le \frac{1}{r}\Longrightarrow & \int_{R_0}^R\frac{\rho'(r)}{\rho(r)}\le \log R-\log R_0 \\ \Longrightarrow & \log(\rho(R))-\log(\rho(R_0))\le \log R-\log R_0 \\ \Longrightarrow & \rho(R)\le CR,\end{split}\end{equation*}

\noindent where $C$ is a constant, expressed in terms of $R_0$. It is a contradiction, because $$\rho(r)\le CR\Longrightarrow \liminf_{r\rightarrow\infty}\rho'(r)\le C$$ and, by hypothesis, $\lim_{r\rightarrow\infty}\rho'(r)=\infty$. Therefore, for $r$ large, $r\rho(r)^{-1}\rho'(r)>1$ and everything else follows in a manner analogous to the proof of \Cref{cit3.3}. The proof of (ii) follows from (i) taking the analysis of $\tilde{\rho}$ defined by $\tilde{\rho}(r)=\rho(-r)$.

\noindent $\blacksquare$

\noindent\textbf{Proof of \Cref{cit3.5}}. Using the \Cref{cit3.1}, we only have to prove that for $r$ of compact support only dependent on $r$ and $b,c$ such that $\pi_1(supp f)\subset [b,c]$:

\vspace{10pt}

\noindent$\displaystyle\int_b^cf_r^2\rho(r)^{n-1}dr+a(n-1)(n-2)\int_b^c(\rho'(r)^2+1)\rho(r)^{n-3}f^2dr$

\noindent$\displaystyle+4a(n-1)\int_b^c\rho'(r)\rho(r)^{n-2}ff_rdr\ge 0$

\noindent$\Longleftrightarrow \displaystyle\int_b^cf_r^2\rho(r)^{n-1}dr+4a(n-1)\int_b^c\rho'(r)\rho(r)^{n-2}ff_rdr+4a^2(n-1)^2\int_b^c\rho'(r)^2\rho(r)^{n-3}f^2dr$

\noindent$\displaystyle-4a^2(n-1)^2\int_b^c\rho'(r)^2\rho(r)^{n-3}f^2dr+a(n-1)(n-2)\int_b^c(\rho'(r)^2+1)\rho(r)^{n-3}f^2dr\ge 0$

\noindent$\displaystyle\Longleftrightarrow \int_b^c(f_r\rho(r)^{\frac{n-1}{2}}+2a(n-1)\rho'(r)\rho(r)^{\frac{n-3}{2}}f)^2dr$

\noindent$\displaystyle+\int_b^c\left(-4a(n-1)\rho'(r)^2+(n-2)(\rho'(r)^2+1)\right)a(n-1)\rho(r)^{n-3}f^2dr\ge 0$.

\vspace{10pt}

\noindent Because the first term of the least inequality is nonnegative, the inequality holds if
$$ -4a(n-1)\rho'(r)^2+(n-2)(\rho'(r)^2+1)\ge 0\Longleftrightarrow a\le \frac{(\rho'(r)^2+1)(n-2)}{4\rho'(r)^2(n-1)}.$$

\noindent Since the application $r\mapsto\frac{1+r}{r}$ is decreasing and $\rho'(r)\le C$ for all $r>0$, we conclude that $M$ is $a$-stable for all $0\le a\le\frac{(C^2+1)(n-2)}{4C^2(n-1)}$.

\noindent$\blacksquare$

\noindent\textbf{Proof of \Cref{cit3.8}}. \hypertarget{proof3.1.5}{} Taking $\rho(r)=r^{\alpha}\xi(r)$, then $C_1\le \xi(r)\le C_2$ for $r\ge 1$. Let $R>0$ be large and the family $f_{R,\beta}$ of functions defined by:

$$f(r)=f_{R,\beta}(r)=\begin{cases}\hspace{40mm}0\hspace{40mm} \hspace{4mm} \text{if} \hspace{4mm} 0\le r\le \frac{m-1}{m}R; \\ \hspace{6.5mm}(2R^{-\frac{n\alpha-\alpha+1}{2}}r-R^{-\frac{n\alpha-\alpha-1}{2}})\xi(R)^{-2a(n-1)} \hspace{7.5mm} \hspace{4mm}\text{if} \hspace{4mm} \frac{R}{2}\le r\le R; \\ \hspace{23mm}r^{-\frac{n\alpha-\alpha-1}{2}}\xi(r)^{-2a(n-1)}\hspace{27mm} \text{if} \hspace{4mm} R\le r\le R^{\beta}; \\ (-2R^{-\frac{n\alpha\beta-\alpha\beta+\beta}{2}}r+3R^{-\frac{n\alpha\beta-\alpha\beta-\beta}{2}})\xi(r)^{-2a(n-1)}\hspace{7mm}\text{if} \hspace{4mm} R^{\beta}\le r\le \frac{3}{2}R^{\beta}; \\ \hspace{40mm}0\hspace{40mm} \hspace{4mm}\text{if}\hspace{4mm} r\ge \frac{3}{2}R^{\beta}.\end{cases}$$

\noindent Let us analyze the value of the expression obtained in (\ref{2.4}) in each one of the intervals $[\frac{R}{2},R]$, $[R,R^{\beta}]$, $[R^{\beta},\frac{3}{2}R^{\beta}]$ separately.

\noindent For interval $[\frac{R}{2},R]$, we will only define the value of (\ref{2.4})  by $K_R$. Let us now analyze the expression (\ref{2.4}) in the interval $[R,R^{\beta}]$. In this interval, we have:

\noindent (i) $\displaystyle f(r)=r^{-\frac{n\alpha-\alpha-1}{2}}\xi(r)^{-2a(n-1)}$;

\noindent (ii) $\displaystyle f_r=-\frac{n\alpha-\alpha-1}{2}r^{-\frac{n\alpha-\alpha+1}{2}}\xi(r)^{-2a(n-1)}-2a(n-1)r^{-\frac{n\alpha-\alpha-1}{2}}\xi(r)^{-2a(n-1)-1}\xi'(r)$;

\noindent (iii) $\displaystyle f_r^2=\frac{(n\alpha-\alpha-1)^2}{4}r^{-n\alpha+\alpha-1}\xi(r)^{-4a(n-1)}$

\hspace{12mm}$+2a(n-1)(n\alpha-\alpha-1)r^{-n\alpha+\alpha}\xi(r)^{-4a(n-1)-1}\xi'(r)$

\hspace{12mm}$+4a^2(n-1)^2r^{-n\alpha+\alpha+1}\xi(r)^{-4a(n-1)-2}\xi'(r)^2$; 

\noindent (iv) $\displaystyle ff_r=-\frac{n\alpha-\alpha-1}{2}r^{-n\alpha+\alpha}\xi(r)^{-4a(n-1)}-2a(n-1)r^{-n\alpha+\alpha+1}\xi(r)^{-4a(n-1)-1}\xi'(r).$ 

\noindent Substituting the expression of $f=f_{R,\beta}$ in (\ref{2.4}) in the interval $[R,R^{\beta}]$ and using (i), (ii), (iii) and (iv), we obtain

\vspace{10pt}

\noindent $\displaystyle\frac{(n\alpha-\alpha-1)^2}{4}\int_R^{R^{\beta}}r^{-1}\xi(r)^{(1-4a)(n-1)}dr+2a(n-1)(n\alpha-\alpha-1)\int_R^{R^{\beta}}\xi(r)^{(1-4a)(n-1)-1}\xi'(r)dr$

\noindent$\displaystyle+4a^2(n-1)^2\int_R^{R^{\beta}}r\xi(r)^{(1-4a)(n-1)-2}\xi'(r)^2dr+a\alpha^2(n-1)(n-2)\int_R^{R^{\beta}}r^{-1}\xi(r)^{(1-4a)(n-1)}dr$

\noindent$\displaystyle +2a\alpha(n-1)(n-2)\int_R^{R^{\beta}}\xi(r)^{(1-4a)(n-1)-1}\xi'(r)dr$

\noindent$\displaystyle+a(n-1)(n-2)\int_R^{R^{\beta}}r\xi(r)^{(1-4a)(n-1)-2}\xi'(r)^2dr$

\noindent$\displaystyle-2a\alpha(n-1)(n\alpha-\alpha-1)\int_R^{R^{\beta}}r^{-1}\xi(r)^{(1-4a)(n-1)}dr$

\noindent$\displaystyle-8a^2\alpha(n-1)^2\int_R^{R^{\beta}}\xi(r)^{(1-4a)(n-1)-1}\xi'(r)dr$

\noindent$\displaystyle-2a(n-1)(n\alpha-\alpha-1)\int_R^{R^{\beta}}\xi(r)^{(1-4a)(n-1)-1}\xi'(r)dr$

\noindent$\displaystyle-8a^2(n-1)^2\int_R^{R^{\beta}}r\xi(r)^{(1-4a)(n-1)-2}\xi'(r)^2dr$

\noindent$\displaystyle+a\frac{S(F)}{A(F)}\int_{R}^{R^{\beta}}r^{-2\alpha+1}\xi(r)^{(1-4a)(n-1)-2}dr$

\noindent$\displaystyle=\left(\frac{(n\alpha-\alpha-1)^2}{4}-a\alpha(n-1)(n\alpha-2)\right)\int_R^{R^{\beta}}r^{-1}\xi(r)^{(1-4a)(n-1)}dr$

\noindent$\displaystyle+2a\alpha(n-1)[(1-4a)(n-1)-1]\int_R^{R^{\beta}}\xi(r)^{(1-4a)(n-1)-1}\xi'(r)dr$

\noindent$\displaystyle+a(n-1)[(1-4a)(n-1)-1]\int_R^{R^{\beta}}r\xi(r)^{(1-4a)(n-1)-2}\xi'(r)^2dr$

\noindent$\displaystyle+a\frac{S(F)}{A(F)}\int_{R}^{R^{\beta}}r^{-2\alpha+1}\xi(r)^{(1-4a)(n-1)-2}dr.$

\vspace{-4mm}

\begin{center}\hypertarget{16}{(16)}\end{center}
Note that, since $a>\frac{n-2}{4(n-1)}$, $$(1-4a)(n-1)-1<\left(1-\frac{n-2}{n-1}\right)(n-1)-1=0,$$ $(1-4a)(n-1)-1$ is negative. By hypothesis, $a>\frac{(n\alpha-\alpha-1)^2}{4\alpha(n-1)(n\alpha-2)}$, therefore, the first and third terms of (\hyperlink{16}{16}) are negatives. The fourth term of (\hyperlink{16}{16}) is significantly small when $R$ is large because $-2\alpha+1<-1$, implying that this term is less than $\frac{S(F)}{(2\alpha-2)A(F)}R^{2-2\alpha}\rightarrow 0$ when $R\rightarrow\infty$. The second term of (\hyperlink{16}{16}) is equal to $$\frac{2a\alpha[(1-4a)(n-1)-1]}{1-4a}(\xi(R^{\beta})^{(1-4a)(n-1)}-\xi(R)^{(1-4a)(n-1)}),$$ if $a\ne\frac{1}{4}$, and equal to $$2a\alpha[(1-4a)(n-1)-1][\log(\xi(R^{\beta}))-\log(\xi(R))],$$ if $a=\frac{1}{4}$.

\noindent In both cases, it is limited because in $[R,R^{\beta}]$, $C_1\le \xi(r)\le C_2$. Since there exists $C_3>0$ such that $\xi(r)^{(1-4a)(n-1)}\ge C_3$ for all $r\in[R,R^{\beta}]$, where we can take $$C_3=\min\{C_1^{(1-4a)(n-1)},C_2^{(1-4a)(n-1)}\},$$ then $$\int_R^{R^{\beta}}r^{-1}\xi(r)^{(1-4a)(n-1)}dr\ge C_3(\beta-1)\log R$$ can be so large as we want. Therefore, the expression in (\ref{2.4}) in the interval $[R,R^{\beta}]$ can be so (negatively) large as we want, based on the choice of $\beta$.

Because of the previous conclusion, on interval $[R^{\beta},\frac{3}{2}R^{\beta}]$, we only have to prove that the expression (\ref{2.4}) in this interval is upper limited for a constant that is independent of $\beta$. We have $$f(r)=(-2R^{-\frac{n\alpha\beta-\alpha\beta+\beta}{2}}r+3R^{-\frac{n\alpha\beta-\alpha\beta-\beta}{2}})\xi(r)^{-2a(n-1)},$$ then $$f_r(r,\cdot)=-2R^{-\frac{n\alpha\beta-\alpha\beta+\beta}{2}}\xi(r)^{-2a(n-1)}-2a(n-1)f\xi(r)^{-1}\xi'(r).$$ Substituting in (\ref{2.4}) and using that there is a constant $C_2>0$ such that $\xi(r)\le C_2$ for all $r\ge R^{\beta}$, we obtain that the expression (\ref{2.4}) in interval $[R^{\beta},\frac{3}{2}R^{\beta}]$, when $\beta\rightarrow\infty$, is:

\vspace{10pt}

\noindent$\displaystyle \underbrace{4R^{-n\alpha\beta+\alpha\beta-\beta}\int_{R^{\beta}}^{R^{\frac{3}{2}\beta}}r^{n\alpha-\alpha}\xi(r)^{(1-4a)(n-1)}dr}_{O(1)}
+4a^2(n-1)^2\int_{R^{\beta}}^{\frac{3}{2}R^{\beta}}r^{n\alpha-\alpha}f^2\xi(r)^{n-3}\xi'(r)^2dr$

\noindent $\displaystyle+8a(n-1)R^{-\frac{n\alpha\beta-\alpha\beta+\beta}{2}}\int_{R^{\beta}}^{\frac{3}{2}R^\beta}r^{n\alpha-\alpha}f\xi(r)^{(1-2a)(n-1)-1}\xi'(r)dr$

\noindent$\displaystyle +\underbrace{4a\alpha^2(n-1)(n-2)R^{-n\alpha\beta+\alpha\beta-\beta}\int_{R^\beta}^{\frac{3}{2}R^{\beta}}r^{n\alpha-\alpha}\xi(r)^{(1-4a)(n-1)}dr}_{O(1)}$

\noindent$\displaystyle-\underbrace{12a\alpha^2(n-1)(n-2)R^{-n\alpha\beta+\alpha\beta}\int_{R^\beta}^{\frac{3}{2}R^{\beta}}r^{n\alpha-\alpha-1}\xi(r)^{(1-4a)(n-1)}dr}_{O(1)}$

\noindent$\displaystyle+\underbrace{9a\alpha^2(n-1)(n-2)R^{-n\alpha\beta+\alpha\beta+\beta}\int_{R^{\beta}}^{\frac{3}{2}R^{\beta}}r^{n\alpha-\alpha-2}\xi(r)^{(1-4a)(n-1)}dr}_{O(1)}$

\noindent$\displaystyle+2a\alpha(n-1)(n-2)\int_{R^{\beta}}^{\frac{3}{2}R^{\beta}}r^{n\alpha-\alpha-1}f^2\xi(r)^{n-2}\xi'(r)dr$

\noindent$\displaystyle+a(n-1)(n-2)\int_{R^{\beta}}^{\frac{3}{2}R^{\beta}}r^{n\alpha-\alpha}f^2\xi(r)^{n-3}\xi'(r)^2dr$

\noindent$\displaystyle-\underbrace{8a\alpha(n-1)R^{-\frac{n\alpha\beta-\alpha\beta+\beta}{2}}\int_{R^{\beta}}^{\frac{3}{2}R^{\beta}}r^{n\alpha-\alpha-1}\xi(r)^{(1-2a)(n-1)}fdr}_{O(1)}$

\noindent$\displaystyle-8a^2\alpha(n-1)^2\int_{R^{\beta}}^{\frac{3}{2}R^{\beta}}r^{n\alpha-\alpha-1}f^2\xi(r)^{n-2}\xi'(r)dr$

\noindent $\displaystyle-8a(n-1)R^{-\frac{n\alpha\beta-\alpha\beta+\beta}{2}}\int_{R^{\beta}}^{\frac{3}{2}R^{\beta}}r^{n\alpha-\alpha}f\xi(r)^{(1-2a)(n-1)-1}\xi'(r)dr$

\noindent$\displaystyle-8a^2(n-1)^2\int_{R^{\beta}}^{\frac{3}{2}R^{\beta}}r^{n\alpha-\alpha}f^2\xi(r)^{n-3}\xi'(r)^2dr+\underbrace{a\frac{S(F)}{A(F)}\int_{R^{\beta}}^{\frac{3}{2}R^{\beta}}r^{n\alpha-3\alpha}\xi(r)^{n-3}f^2dr}_{o(1)}$

\noindent$\displaystyle\le C+a(n-1)\underbrace{(n-2-4(n-1))}_{\le 0}\int_{R^{\beta}}^{\frac{3}{2}R^{\beta}}r^{n\alpha-\alpha}f^2\xi(r)^{n-3}\xi'(r)^2dr$

\noindent$\displaystyle+2a\alpha(n-1)(n-2-4a(n-1))\int_{R^{\beta}}^{\frac{3}{2}R^{\beta}}r^{n\alpha-\alpha-1}f^2\xi(r)^{n-2}\xi'(r)dr$.

\vspace{10pt}

\noindent It remains to be shown that the last term is limited by a constant independent of $\beta$. If $a\ne\frac{1}{4}$, we have:
\begin{equation*}\begin{split}\int_{R^{\beta}}^{\frac{3}{2}R^{\beta}}r^{n\alpha-\alpha-1}f^2\xi(r)^{n-2}\xi'(r)dr= & 4R^{-n\alpha\beta+\alpha\beta-\beta}\int_{R^{\beta}}^{\frac{3}{2}R^{\beta}}r^{n\alpha-\alpha+1}\xi(r)^{(1-4a)(n-1)-1}\xi'(r)dr \\ & -12R^{-n\alpha\beta+\alpha\beta}\int_{R^{\beta}}^{\frac{3}{2}R^{\beta}}r^{n\alpha-\alpha}\xi(r)^{(1-4a)(n-1)-1}\xi'(r)dr \\ & +9R^{-n\alpha\beta+\alpha\beta+\beta}\int_{R^{\beta}}^{\frac{3}{2}R^{\beta}}r^{n\alpha-\alpha-1}\xi(r)^{(1-4a)(n-1)-1}\xi'(r)dr \\ &=
\frac{4R^{-n\alpha\beta+\alpha\beta-\beta}}{(1-4a)(n-1)}\left[ r^{n\alpha-\alpha+1}\xi(r)^{(1-4a)(n-1)}\left.\right|_{R^{\beta}}^{\frac{3}{2}R^{\beta}}\phantom{\int_{R^{\beta}}^{\frac{3}{2}R^{\beta}}}\right.  \\ & \left.-(n\alpha-\alpha+1)\int_{R^{\beta}}^{\frac{3}{2}R^{\beta}}r^{n\alpha-\alpha}\xi(r)^{(1-4a)(n-1)}dr\right] \\  & -\frac{12R^{-n\alpha\beta+\alpha\beta}}{(1-4a)(n-1)}\left[ r^{n\alpha-\alpha}\xi(r)^{(1-4a)(n-1)}\left.\right|_{R^{\beta}}^{\frac{3}{2}R^{\beta}}\phantom{\int_{R^{\beta}}^{\frac{3}{2}R^{\beta}}} \right. \end{split}\end{equation*}
\begin{equation*}\begin{split}  \phantom{\int_{R^{\beta}}^{\frac{3}{2}R^{\beta}}r^{n\alpha-\alpha-1}f^2\xi(r)^{n-2}\xi'(r)dr= } & \left.-(n\alpha-\alpha)\int_{R^{\beta}}^{\frac{3}{2}R^{\beta}}r^{n\alpha-\alpha-1}\xi(r)^{(1-4a)(n-1)}dr\right] \\ & +\frac{9R^{-n\alpha\beta+\alpha\beta+\beta}}{(1-4a)(n-1)}\left[ r^{n\alpha-\alpha-1}\xi(r)^{(1-4a)(n-1)}\left.\right|_{R^{\beta}}^{\frac{3}{2}R^{\beta}}\phantom{\int_{R^{\beta}}^{\frac{3}{2}R^{\beta}}}\right. \\ & \left.-(n\alpha-\alpha-1)\int_{R^{\beta}}^{\frac{3}{2}R^{\beta}}r^{n\alpha-\alpha-2}\xi(r)^{(1-4a)(n-1)}dr\right].\end{split}\end{equation*}
\noindent Using that $C_1\le\xi(r)\le C_2$, we can conclude that the last expression is limited by a constant independent of $\beta$. This ends the proof for the case $a\ne\frac{1}{4}$. If $a=\frac{1}{4}$, we have

\vspace{10pt}

\noindent\hspace{13.4mm}$\displaystyle \int_{R^{\beta}}^{\frac{3}{2}R^{\beta}}r^{n\alpha-\alpha-1}f^2\xi(r)^{n-2}\xi'(r)dr=4R^{-n\alpha\beta+\alpha\beta-\beta}\int_{R^{\beta}}^{\frac{3}{2}R^{\beta}}r^{n\alpha-\alpha+1}\xi(r)^{-1}\xi'(r)dr$\hspace{13.4mm}

\hspace{13.4mm}\hspace{54mm}$\displaystyle -12R^{-n\alpha\beta+\alpha\beta}\int_{R^{\beta}}^{\frac{3}{2}R^{\beta}}r^{n\alpha-\alpha}\xi(r)^{-1}\xi'(r)dr$

\hspace{13.4mm}\hspace{54mm}$\displaystyle+9R^{-n\alpha\beta+\alpha\beta+\beta}\int_{R^{\beta}}^{\frac{3}{2}R^{\beta}}r^{n\alpha-\alpha-1}\xi(r)^{-1}\xi'(r)dr$

\hspace{13.4mm}\hspace{50mm}$\displaystyle=4R^{-n\alpha\beta+\alpha\beta-\beta}\left[ r^{n\alpha-\alpha-1}\log(\xi(r))\left.\right|_{R^{\beta}}^{\frac{3}{2}R^{\beta}}\phantom{\int_{R^{\beta}}^{\frac{3}{2}R^{\beta}}}\right.$

\hspace{13.4mm}\hspace{54mm}$\displaystyle\left.-(n\alpha-\alpha+1)\int_{R^{\beta}}^{\frac{3}{2}R^{\beta}}r^{n\alpha-\alpha}\log(\xi(r))dr\right]$

\hspace{13.4mm}\hspace{54mm}$\displaystyle-12R^{-n\alpha\beta+\alpha\beta}\left[ r^{n\alpha-\alpha}\log(\xi(r))\left.\right|_{R^{\beta}}^{\frac{3}{2}R^{\beta}}\phantom{\int_{R^{\beta}}^{\frac{3}{2}R^{\beta}}}\right.$

\hspace{13.4mm}\hspace{54mm}$\displaystyle\left.-(n\alpha-\alpha)\int_{R^{\beta}}^{\frac{3}{2}R^{\beta}}r^{n\alpha-\alpha-1}\log(\xi(r))dr\right]$

\hspace{13.4mm}\hspace{54mm}$\displaystyle+9R^{-n\alpha\beta+\alpha\beta+\beta}\left[ r^{n\alpha-\alpha-1}\log(\xi(r))\left.\right|_{R^{\beta}}^{\frac{3}{2}R^{\beta}}\phantom{\int_{R^{\beta}}^{\frac{3}{2}R^{\beta}}}\right.$

\hspace{13.4mm}\hspace{54mm}$\displaystyle\left.-(n\alpha-\alpha-1)\int_{R^{\beta}}^{\frac{3}{2}R^{\beta}}r^{n\alpha-\alpha-2}\log(\xi(r))dr\right]$.

\vspace{10pt}

\noindent As in the previous case, the limitation $C_1\le \xi(r)\le C_2$ ensures that the expression is limited by a constant independent of $\beta$. This completes the proof. 

In the next proposition, we show that the constant $\frac{(n\alpha-\alpha-1)^2}{4\alpha(n-1)(n\alpha-2)}$ which lowers the constant $a$ is the best possible in the general case of \Cref{cit3.8}.

\begin{proposicao}\label{cit3.23}
    Let $n\ge 2$ and $M$ be the $n$-dimensional plane with the metric $dr^2+\rho(r)^2g_{\mathbb S^{n-1}}$, when $\rho(r)=\frac{1}{\alpha+1}r^{\alpha}$, with $\alpha>1$, then $L$ is $a$-stable for all $a\le \frac{(n\alpha-\alpha-1)^2}{4\alpha(n-1)(n\alpha-2)}$. 
\end{proposicao}

\noindent\textbf{Proof}. Note that, for all real function $f$ of compact support:
\begin{equation*}\begin{split}\int_0^{\infty}r^{n\alpha-\alpha-2}f^2dr= & -\frac{2}{n\alpha-\alpha-1}\int_0^{\infty}r^{n\alpha-\alpha-1}ff'dr\\ = & -\frac{2}{n\alpha-\alpha-1}\int_0^{\infty}r^{\frac{n\alpha-\alpha-2}{2}}fr^{\frac{n\alpha-\alpha}{2}}f'dr\\ \le & \frac{2}{n\alpha-\alpha-1}\left(\int_0^{\infty}r^{n\alpha-\alpha-2}f^2dr\right)^{\frac{1}{2}}\left(\int_0^{\infty}r^{n\alpha-\alpha}(f')^2dr\right)^{\frac{1}{2}}.\end{split}\end{equation*}

\noindent Therefore
$$\int_0^{\infty}r^{n\alpha-\alpha-2}f^2dr\le \frac{4}{(n\alpha-\alpha-1)^2}\int_0^{\infty}r^{n\alpha-\alpha}(f')^2dr.$$

\noindent By \Cref{cit3.1}, we can reduce the space of test functions to the space of smooth compact support functions that depend only on $r$. Suppose $f$ is dependent only on $r$, 
substituting in (\hyperlink{3.8}{3.8}) (note that, in this case, $\xi(r)=\frac{1}{\alpha+1}$), we have

\vspace{10pt}

\noindent \hspace{1mm}$\displaystyle\frac{1}{(\alpha+1)^{n-1}A(F)}\int_M-fL_af=\int_0^{\infty}f_r^2r^{n\alpha-\alpha}dr-a\alpha(n-1)(n\alpha-2)\int_0^{\infty}r^{n\alpha-\alpha-2}f^2dr$

\hspace{1mm}\hspace{49mm}$\displaystyle +a(n-1)(n-2)(\alpha+1)^2\int_0^{\infty}r^{n\alpha-3\alpha}f^2dr$

\hspace{1mm}\hspace{44mm}$\displaystyle\ge \left(\frac{(n\alpha-\alpha-1)^2}{4}-a\alpha(n-1)(n\alpha-2)\right)\int_0^{\infty}r^{n\alpha-\alpha-2}f^2dr$,

\vspace{10pt}

\noindent the last term is nonnegative for all $a\le\frac{(n\alpha-\alpha-1)^2}{4\alpha(n-1)(n\alpha-2)}$.

\noindent$\blacksquare$

\noindent\textbf{Proof of \Cref{cit3.11}}.\hypertarget{proof3.1.6}{} Since $\xi(r)\rightarrow 0$ when $r\rightarrow\infty$ and $r\xi(r)^{-1}\xi'(r)\ge -(\alpha-1)$, we have \begin{equation*}\begin{split} \xi(r)^{-1}\xi'(r) & \ge-\frac{\alpha-1}{r} \\ \Longrightarrow \int_1^r\xi(s)^{-1}\xi'(s)ds & \ge -(\alpha-1)\int_1^rs^{-1}ds\\ \Longrightarrow \log(\xi(r))-\log(\xi(1)) & \ge -(\alpha-1)\log r\\ \Longrightarrow \xi(r) & \ge C_1r^{1-\alpha},\end{split}\end{equation*} where $C_1=\xi(1)^{-1}$. In order to prove by contradiction, let us suppose that $$\limsup_{r\rightarrow\infty}\xi(r)r^{\alpha-1}=\infty.$$ Since $$(\xi(r)r^{\alpha-1})'=\xi'(r)r^{\alpha-1}+(\alpha-1)\xi(r)r^{\alpha-2}\ge -(\alpha-1)r^{-1}\xi(r)r^{\alpha-1}+(\alpha-1)\xi(r)r^{\alpha-2}=0,$$ the application $r\mapsto \xi(r)r^{\alpha-1}$ is nondecreasing and we can assume that $$\lim_{r\rightarrow\infty}\xi(r)r^{\alpha-1}=\infty.$$
Take as a test function
\begin{equation*}f(r)=f_{R,\beta}(r,\cdot)=\left\{\begin{array}{cl}0 & \text{if} \hspace{4mm} 0\le r\le \frac{R}{2}; \\ (2R^{-\frac{n\alpha-\alpha+1}{2}}r-R^{-\frac{n\alpha-\alpha-1}{2}})\xi(R)^{-\frac{1}{2}(n-1)} & \text{if} \hspace{4mm} \frac{R}{2}\le r\le R; \\ r^{-\frac{n\alpha-\alpha-1}{2}}\xi(r)^{-\frac{1}{2}(n-1)} & \text{if} \hspace{4mm} R\le r\le R^{\beta}; \\ (-2R^{-\frac{n\alpha\beta-\alpha\beta+\beta}{2}}r+3R^{-\frac{n\alpha\beta-\alpha\beta-\beta}{2}})\xi(r)^{-\frac{1}{2}(n-1)} & \text{if} \hspace{4mm} R^{\beta}\le r\le \frac{3}{2}R^{\beta}; \\ 0 & \text{if}\hspace{4mm} r\ge \frac{3}{2}R^{\beta}.\end{array}\right.\end{equation*}

\noindent As in the proof of the previous theorem, let us analyze the value of the expression obtained in (\ref{2.4}) in each of the intervals: $[\frac{R}{2},R]$, $[R,R^{\beta}]$, $[R^{\beta}$ and $\frac{3}{2}R^{\beta}]$ separately.

On interval $[\frac{R}{2},R]$, we will only define the value of (\ref{2.4}) on interval $[\frac{R}{2},R]$ by $K_R$. Let us now analyze the expression (\ref{2.4}) in the interval $[R,R^{\beta}]$. In this interval, we have 

\noindent (i) $\displaystyle f(r)=r^{-\frac{n\alpha-\alpha-1}{2}}\xi(r)^{-\frac{1}{2}(n-1)}$;

\noindent (ii) $\displaystyle f_r=-\frac{n\alpha-\alpha-1}{2}r^{-\frac{n\alpha-\alpha+1}{2}}\xi(r)^{-\frac{1}{2}(n-1)}-\frac{1}{2}(n-1)r^{-\frac{n\alpha-\alpha-1}{2}}\xi(r)^{-\frac{1}{2}(n-1)-1}\xi'(r)$;

\noindent (iii) $\displaystyle f_r^2=\frac{(n\alpha-\alpha-1)^2}{4}r^{-n\alpha+\alpha-1}\xi(r)^{-n+1}+\frac{1}{2}(n-1)(n\alpha-\alpha-1)r^{-n\alpha+\alpha}\xi(r)^{-n}\xi'(r)$

\vspace{3mm}\hspace{11mm}$\displaystyle +\frac{1}{4}(n-1)^2r^{-n\alpha+\alpha+1}\xi(r)^{-n-1}\xi'(r)^2$;

\noindent (iv) $\displaystyle ff_r=-\frac{n\alpha-\alpha-1}{2}r^{-n\alpha+\alpha}\xi(r)^{-n+1}-\frac{1}{2}(n-1)r^{-n\alpha+\alpha+1}\xi(r)^{-n}\xi'(r)$. 

\noindent Substituting the expression of $f=f_{R,\beta}$ in (\ref{2.4}) in the interval $[R,R^{\beta}]$ and using (i), (ii), (iii) and (iv), we obtain

\vspace{10pt}

\noindent $\displaystyle\frac{(n\alpha-\alpha-1)^2}{4}\int_R^{R^{\beta}}r^{-1}dr+\frac{1}{2}(n-1)(n\alpha-\alpha-1)\int_R^{R^{\beta}}\xi(r)^{-1}\xi'(r)dr$

\noindent$\displaystyle+\frac{1}{4}(n-1)^2\int_R^{R^{\beta}}r\xi(r)^{-2}\xi'(r)^2dr+a\alpha^2(n-1)(n-2)\int_R^{R^{\beta}}r^{-1}dr$

\noindent$\displaystyle +2a\alpha(n-1)(n-2)\int_R^{R^{\beta}}\xi(r)^{-1}\xi'(r)dr+a(n-1)(n-2)\int_R^{R^{\beta}}r\xi(r)^{-2}\xi'(r)^2dr$

\noindent$\displaystyle-2a\alpha(n-1)(n\alpha-\alpha-1)\int_R^{R^{\beta}}r^{-1}dr-2a\alpha(n-1)^2\int_R^{R^{\beta}}\xi(r)^{-1}\xi'(r)dr$

\noindent$\displaystyle-2a(n-1)(n\alpha-\alpha-1)\int_R^{R^{\beta}}\xi(r)^{-1}\xi'(r)dr-2a(n-1)^2\int_R^{R^{\beta}}r\xi(r)^{-2}\xi'(r)^2dr$

\noindent$\displaystyle+a\frac{S(F)}{A(F)}\int_{R}^{R^{\beta}}r^{-2\alpha+1}\xi(r)^{-2}dr$

\noindent$\displaystyle=\left(\frac{(n\alpha-\alpha-1)^2}{4}-a\alpha(n-1)(n\alpha-2)\right)\int_R^{R^{\beta}}r^{-1}dr$

\noindent$\displaystyle+\frac{1}{2}(n-1)(n\alpha-\alpha-4a\alpha n+4a-1)\int_R^{R^{\beta}}\xi(r)^{-1}\xi'(r)dr$

\noindent$\displaystyle+(n-1)\left(\frac{1}{4}n-an-\frac{1}{4}\right)\int_R^{R^{\beta}}r\xi(r)^{-2}\xi'(r)^2dr$

\noindent$\displaystyle+a\frac{S(F)}{A(F)}\int_{R}^{R^{\beta}}r^{-2\alpha+1}\xi(r)^{-2}dr.\hspace{109.4mm}\hypertarget{17}{(17)}$

\vspace{10pt}

\noindent By hypothesis, $a\ge\frac{(n\alpha-\alpha-1)^2}{4\alpha(n-1)(n\alpha-2)}$, then the first term of (\hyperlink{17}{17}) is nonpositive. 
Let us show that the second term of (\hyperlink{17}{17}) increases arbitrarily more than the fourth term of (\hyperlink{17}{17}) (in module), if we choose $R$ and $\beta$ appropriately. 
Note that $$\int_R^{R^{\beta}}\xi(r)^{-1}\xi'(r)dr=\log(\xi(R^{\beta}))-\log(\xi(r))\rightarrow\infty$$
when $R\rightarrow\infty$, because $\xi(r)\rightarrow0$ when $r\rightarrow\infty$. If $$\int_R^{\infty}r^{-2\alpha+1}\xi(r)^{-2}dr<\infty,$$ there is nothing to do. Suppose that $$\int_R^{\infty}r^{-2\alpha+1}\xi(r)^{-2}dr=\infty,$$ then for all $\delta>0$, $$\liminf_{r\rightarrow\infty}\xi(r)r^{\alpha-1-\delta}=0,$$ because otherwise, $$\liminf_{r\rightarrow\infty}\xi(r)r^{\alpha-1-\delta}=C>0,$$ where $C$ can be $\infty$. It implies $$\limsup_{r\rightarrow\infty}\xi(r)^{-2}r^{-2\alpha+2+2\delta}= C^{-2}.$$ Then for $R$ sufficiently large, $$\int_R^{\infty}r^{-2\alpha+1}\xi(r)^{-2}dr=\int_R^{\infty}r^{-1-2\delta}\xi(r)^{-2}r^{-2\alpha+2+2\delta}\le (C^{-2}+1)\int_R^{\infty}r^{-1-2\delta}<\infty.$$ It is a contradiction. Hence, there exists a sequence $\{R_1,R_2,\dots\}$ such that $$\xi(R_i)<R_i^{-\alpha+1+\frac{\alpha-1}{2}}\Longrightarrow \xi(R_i)<R_i^{-\frac{\alpha-1}{2}}\Longrightarrow \log(\xi(R_i))<-\frac{\alpha-1}{2}\log R_i$$ for all positive integer $i$. Therefore, it is possible to choose $\beta$ large and appropriately such that $$\int_R^{R^{\beta}}\xi(r)^{-1}\xi'(r)dr<C_1\log (R^{\beta})=C_1\beta\log R$$ for some negative constant $C_1$ (we can choose $\beta$ such that $R^{\beta}\in \{R_1,R_2,\dots\}$, where $\{R_i\}$ is the sequence defined above and $C_1$ a little bigger  than $-\frac{\alpha-1}{2}$). On the other hand,$$\int_r^{R^{\beta}}r^{-2\alpha+1}\xi(r)^{-2}dr=\int_R^{R^{\beta}}(\xi(r)^{-2}r^{-2\alpha+2})r^{-1}dr.$$ Since we are assuming that $\lim_{r\rightarrow\infty}\xi(r)r^{\alpha-1}=\infty$, taking $R$ large, then $$\int_R^{R^{\beta}}r^{-2\alpha+1}\xi(r)^{-2}dr\ll\beta\log R.$$ It shows that  the fourth term of (\hyperlink{17}{17}) is arbitrarily less than the second term of (\hyperlink{17}{17}) (in module). Intending to compare the second and the third term, we have
$$\int_R^{R^{\beta}}r\xi(r)^{-2}\xi'(r)^2dr\le -(\alpha-1)\int_R^{R^{\beta}}\xi(r)^{-1}\xi'(r)dr,$$ then the combination between the second and the third term of (\hyperlink{17}{17}) is less than $$\int_R^{R^{\beta}}\xi(r)^{-1}\xi'(r)dr$$ times the constant 
$$\frac{1}{2}(n-1)(n\alpha-\alpha-4a\alpha n+4a-1)-(\alpha-1)(n-1)\left(\frac{1}{4}n-an-\frac{1}{4}\right)$$ \begin{equation}\label{120} =\frac{1}{4}(n-1)(n\alpha-\alpha-4a\alpha n+8a+n-4an-3).\tag{18}\end{equation} Substituting $a=\frac{(n\alpha-\alpha-1)^2}{4\alpha(n-1)(n\alpha-2)}$ on (\ref{120}), we have

\vspace{10pt}

\noindent $\displaystyle\frac{1}{4}(n-1)\left[\frac{n\alpha^2(n-1)(n\alpha-2)}{\alpha(n-1)(n\alpha-2)}-\frac{\alpha^2(n-1)(n\alpha-2)}{\alpha(n-1)(n\alpha-2)}-\frac{\alpha n(n\alpha-\alpha-1)^2}{\alpha(n-1)(n\alpha-2)}\right.$

\hspace{9mm}
\noindent$\displaystyle+\frac{2(n\alpha-\alpha-1)^2}{\alpha(n-1)(n\alpha-2)}+\frac{n\alpha (n-1)(n\alpha-2)}{\alpha(n-1)(n\alpha-2)}-\frac{n(n\alpha-\alpha-1)^2}{\alpha(n-1)(n\alpha-2)}$

\hspace{9mm}\noindent$\displaystyle\left.-\frac{3\alpha(n-1)(n\alpha-2)}{\alpha(n-1)(n\alpha-2)}\right]$

\noindent$\displaystyle=\frac{1}{4\alpha(n\alpha-2)}(n^3\alpha^3-2n^2\alpha^2-n^2\alpha^3+2n\alpha^2-n^2\alpha^3+2n\alpha^2+n\alpha^3-2\alpha^2-n^3\alpha^3-n\alpha^3-n\alpha$

\hspace{20mm}\noindent$+2n^2\alpha^3+2n^2\alpha^2-2n\alpha^2+2n^2\alpha^2+2\alpha^2+2-4n\alpha^2-4n\alpha+4\alpha+n^3\alpha^2-2n^2\alpha$

\hspace{20mm}\noindent$-n^2\alpha^2+2n\alpha-n^3\alpha^2-n\alpha^2-n+2n^2\alpha^2+2n^2\alpha-2n\alpha-3n^2\alpha^2+6n\alpha+3n\alpha^2-6\alpha)$

\noindent $\displaystyle=\frac{n\alpha-2\alpha-n+2}{4\alpha(n\alpha-2)}=\frac{(n-2)(\alpha-1)}{4\alpha(n\alpha-2)}>0$;

\vspace{10pt}

\noindent because $\alpha>1$ and $n\ge 3$. Therefore, the combination of the second and the third term of (\hyperlink{17}{17}) is less than $$\frac{(n-2)(\alpha-1)}{4\alpha(n\alpha-2)}\int_R^{R^{\beta}}\xi(r)^{-1}\xi'(r)dr.$$ Since the first term of (\hyperlink{17}{17}) is zero and the fourth term is irrelevant to the second term for $R$ large, then, for $R$ sufficiently large, the expression (\hyperlink{17}{17}) is less than $$\frac{(n-2)(\alpha-1)}{8\alpha(n\alpha-2)}\int_R^{R^{\beta}}\xi(r)^{-1}\xi'(r)dr=\frac{(n-2)(\alpha-1)}{8\alpha(n\alpha-2)}(\log(\xi(R^{\beta})-\xi(R))).$$
\noindent Using that $\lim_{r\rightarrow\infty}\xi(r)=0$, we conclude that (\hyperlink{17}{17}) can be (negatively) so large as we want, taking $\beta$ large.

Based on the previous conclusion, on interval $[R^{\beta},\frac{3}{2}R^{\beta}]$, we only have to prove that the substitution of (\ref{2.4}) by $f$ is limited and independent of $\beta$. We have $$f(r)=(-2R^{-\frac{n\alpha\beta-\alpha\beta+\beta}{2}}r+3R^{-\frac{n\alpha\beta-\alpha\beta-\beta}{2}})\xi(r)^{-\frac{1}{2}(n-1)},$$ then $$f_r=-2R^{-\frac{n\alpha\beta-\alpha\beta+\beta}{2}}\xi(r)^{-\frac{1}{2}(n-1)}-\frac{1}{2}(n-1)f\xi(r)^{-1}\xi'(r).$$ Substituting in (\ref{2.4}), we obtain when $\beta\rightarrow\infty$:

\vspace{10pt}

\noindent$\displaystyle \underbrace{4R^{-n\alpha\beta+\alpha\beta-\beta}\int_{R^{\beta}}^{R^{\frac{3}{2}\beta}}r^{n\alpha-\alpha}dr}_{O(1)}
+\frac{1}{4}(n-1)^2\int_{R^{\beta}}^{\frac{3}{2}R^{\beta}}r^{n\alpha-\alpha}f^2\xi(r)^{n-3}\xi'(r)^2dr$

\noindent $\displaystyle+2(n-1)R^{-\frac{n\alpha\beta-\alpha\beta+\beta}{2}}\int_{R^{\beta}}^{\frac{3}{2}R^\beta}r^{n\alpha-\alpha}f\xi(r)^{\frac{1}{2}(n-1)-1}\xi'(r)dr$

\noindent$\displaystyle +\underbrace{4a\alpha^2(n-1)(n-2)R^{-n\alpha\beta+\alpha\beta-\beta}\int_{R^\beta}^{\frac{3}{2}R^{\beta}}r^{n\alpha-\alpha}dr}_{O(1)}$

\noindent$\displaystyle-\underbrace{12a\alpha^2(n-1)(n-2)R^{-n\alpha\beta+\alpha\beta}\int_{R^\beta}^{\frac{3}{2}R^{\beta}}r^{n\alpha-\alpha-1}dr}_{O(1)}$

\noindent$\displaystyle+\underbrace{9a\alpha^2(n-1)(n-2)R^{-n\alpha\beta+\alpha\beta+\beta}\int_{R^{\beta}}^{\frac{3}{2}R^{\beta}}r^{n\alpha-\alpha-2}dr}_{O(1)}$

\noindent$\displaystyle+2a\alpha(n-1)(n-2)\int_{R^{\beta}}^{\frac{3}{2}R^{\beta}}r^{n\alpha-\alpha-1}f^2\xi(r)^{n-2}\xi'(r)dr$

\noindent$\displaystyle+a(n-1)(n-2)\int_{R^{\beta}}^{\frac{3}{2}R^{\beta}}r^{n\alpha-\alpha}f^2\xi(r)^{n-3}\xi'(r)^2dr$

\noindent$\displaystyle-\underbrace{8a\alpha(n-1)R^{-\frac{n\alpha\beta-\alpha\beta+\beta}{2}}\int_{R^{\beta}}^{\frac{3}{2}R^{\beta}}r^{n\alpha-\alpha-1}\xi(r)^{\frac{1}{2}(n-1)}fdr}_{O(1)}$

\noindent$\displaystyle-2a\alpha(n-1)^2\int_{R^{\beta}}^{\frac{3}{2}R^{\beta}}r^{n\alpha-\alpha-1}f^2\xi(r)^{n-2}\xi'(r)dr$

\noindent $\displaystyle-8a(n-1)R^{-\frac{n\alpha\beta-\alpha\beta+\beta}{2}}\int_{R^{\beta}}^{\frac{3}{2}R^{\beta}}r^{n\alpha-\alpha}f\xi(r)^{\frac{1}{2}(n-1)-1}\xi'(r)dr$

\noindent$\displaystyle-2a(n-1)^2\int_{R^{\beta}}^{\frac{3}{2}R^{\beta}}r^{n\alpha-\alpha}f^2\xi(r)^{n-3}\xi'(r)^2dr+\underbrace{a\frac{S(F)}{A(F)}\int_{R^{\beta}}^{\frac{3}{2}R^{\beta}}r^{n\alpha-3\alpha}\xi(r)^{n-3}f^2dr}_{O(1)}$

\noindent$\displaystyle\le C+(n-1)\left(\frac{1}{4}n-an-\frac{1}{4}\right)\int_{R^{\beta}}^{\frac{3}{2}R^{\beta}}r^{n\alpha-\alpha}f^2\xi(r)^{n-3}\xi'(r)^2dr$

\noindent $\displaystyle-2a\alpha(n-1)\int_{R^{\beta}}^{\frac{3}{2}R^{\beta}}r^{n\alpha-\alpha-1}f^2\xi(r)^{n-2}\xi'(r)dr$

\noindent$\displaystyle+(2-8\alpha)(n-1)R^{-\frac{n\alpha\beta-\alpha\beta+\beta}{2}}\int_{R^{\beta}}^{\frac{3}{2}R^{\beta}}r^{n\alpha-\alpha}f\xi(r)^{\frac{1}{2}(n-1)-1}\xi'(r)dr$.\hspace{55.9mm}\hypertarget{2.15}{(19)}

\vspace{10pt}

\noindent We must show that all these terms is upper limited by a constant independent of $\beta$. We can note that, on the interval $[R^{\beta},\frac{3}{2}R^{\beta}]$, $$0\le f(r)\le R^{-\frac{n\alpha\beta-\alpha\beta-\beta}{2}}\xi(r)^{-\frac{1}{2}(n-1)}.$$

\noindent Furthermore, $-(\alpha+1)\le r\xi(r)^{-1}\xi'(r)\le 0$. Hence, \begin{equation*}\begin{split}0\le\int_{R^{\beta}}^{\frac{3}{2}R^{\beta}}r^{n\alpha-\alpha}f^2\xi(r)^{n-3}\xi'(r)^2dr\le & \int_{R^{\beta}}^{\frac{3}{2}R^{\beta}}\left(\frac{r}{R^{\beta}}\right)^{n\alpha-\alpha-1}r\xi(r)^{n-2}\xi'(r)^2dr \\ \le & \left(\frac{3}{2}\right)^{n\alpha-\alpha-1}(\alpha-1)^2\int_{R^{\beta}}^{\frac{3}{2}R^{\beta}}r^{-1}dr \\ =& \left(\frac{3}{2}\right)^{n\alpha-\alpha-1}(\alpha-1)^2\log\left(\frac{3}{2}\right);\end{split}\end{equation*}
\begin{equation*}\begin{split}0\ge\int_{R^{\beta}}^{\frac{3}{2}R^{\beta}}r^{n\alpha-\alpha-1}f^2\xi(r)^{n-2}\xi'(r)dr\ge & \int_{R^{\beta}}^{\frac{3}{2}R^{\beta}}\left(\frac{r}{R^{\beta}}\right)^{n\alpha-\alpha-1}\xi(r)^{-1}\xi'(r)dr \\ \ge & -\left(\frac{3}{2}\right)^{n\alpha-\alpha-1}(\alpha-1)\int_{R^{\beta}}^{\frac{3}{2}R^{\beta}}r^{-1}dr \\ =& -\left(\frac{3}{2}\right)^{n\alpha-\alpha-1}(\alpha-1)\log\left(\frac{3}{2}\right);\end{split}\end{equation*}
\begin{equation*}\begin{split}0\ge R^{-\frac{n\alpha\beta-\alpha\beta+\beta}{2}}\int_{R^{\beta}}^{\frac{3}{2}R^{\beta}}r^{n\alpha-\alpha-1}f\xi(r)^{\frac{1}{2}(n-1)-1}\xi'(r)dr\ge & \int_{R^{\beta}}^{\frac{3}{2}R^{\beta}}\left(\frac{r}{R^{\beta}}\right)^{n\alpha-\alpha-1}\xi(r)^{-1}\xi'(r)dr \\ \ge & -\left(\frac{3}{2}\right)^{n\alpha-\alpha-1}(\alpha-1)\int_{R^{\beta}}^{\frac{3}{2}R^{\beta}}r^{-1}dr \\ =& -\left(\frac{3}{2}\right)^{n\alpha-\alpha-1}(\alpha-1)\log\left(\frac{3}{2}\right).\end{split}\end{equation*}
It shows that (\hyperlink{2.15}{19}) is upper limited by a constant independent of $\beta$ and shows that $M$ is $a$-unstable for $a= \frac{(n\alpha-\alpha-1)^2}{4\alpha(n-1)(n\alpha-2)}$. 

Therefore, $r\mapsto \xi(r)r^{\alpha-1}\le C$ for some $C>0$ and $$\rho(r)=r^{\alpha}\xi(r)\le Cr.$$

\noindent $\blacksquare$

\noindent \textbf{Proof of \Cref{xiiso(1)}} Let us remember (\ref{3.6}), that \begin{equation}\label{e1} \alpha=\inf\{\gamma;\lim_{r\rightarrow\infty}\rho(r)r^{-\gamma}=0\}\tag{20}\end{equation} and $\xi(r)=\rho(r)r^{-\alpha}$. We affirm that $$\limsup_{r\rightarrow\infty}\frac{\log(\xi(r))}{\log r}=0.$$
To verify this, let 
$$C=\limsup_{r\rightarrow\infty}\frac{\log(\xi(r))}{\log r}.$$ If $C>0$, then there exists an unlimited sequence $\{R_1,R_2,\dots\}$ such that $$\frac{\log(\xi(R_i))}{\log R_i}>\frac{C}{2}\Longrightarrow \xi(R_i)>R_i^{\frac{C}{2}} \hspace{3mm}\forall i\in\mathbb N.$$
It contradicts the characterization of $\alpha$ in (\ref{e1}), because it would imply $$\limsup_{r\rightarrow\infty}\rho(r)r^{-(\alpha+\frac{C}{4})}=\limsup_{r\rightarrow\infty}\xi(r)r^{-\frac{C}{4}}\ge\limsup_{i\rightarrow\infty}\xi(R_i)R_i^{-\frac{C}{4}}=+\infty.$$ On the other hand, if $C<0$, there exists $R_0$ such that for all $r>R_0$, $$\frac{\log(\xi(r))}{\log r}<\frac{C}{2}\Longrightarrow \xi(r)<r^{\frac{C}{2}}.$$ It also contradicts the characterization of $\alpha$ in $(\ref{e1})$, because it would imply $$\lim_{r\rightarrow\infty}\rho(r)r^{-(\alpha+\frac{C}{4})}=\lim_{r\rightarrow\infty}\xi(r)r^{-\frac{C}{4}}<\lim_{r\rightarrow\infty}r^{\frac{C}{4}}=0.$$ Note that $\alpha+\frac{C}{4}<\alpha$, contracting the characterization of $\alpha$ in (\ref{e1}). Therefore, fixed $R$, we can choose $\beta$ appropriately large such that \begin{equation}\label{appr}\frac{\int_R^{R^{\beta}}\xi(r)^{-1}\xi'(r)dr}{\beta\log(R)}\approx 0.\tag{21}\end{equation}
Our test function will be the same of the \hyperlink{proof3.1.6}{proof of Theorem 3.1.6}, that is 
\begin{equation*}f(r)=f_{R,\beta}(r,\cdot)=\left\{\begin{array}{cl}0 & \text{if} \hspace{4mm} 0\le r\le \frac{R}{2}; \\ (2R^{-\frac{n\alpha-\alpha+1}{2}}r-R^{-\frac{n\alpha-\alpha-1}{2}})\xi(R)^{-\frac{1}{2}(n-1)} & \text{if} \hspace{4mm} \frac{R}{2}\le r\le R; \\ r^{-\frac{n\alpha-\alpha-1}{2}}\xi(r)^{-\frac{1}{2}(n-1)} & \text{if} \hspace{4mm} R\le r\le R^{\beta}; \\ (-2R^{-\frac{n\alpha\beta-\alpha\beta+\beta}{2}}r+3R^{-\frac{n\alpha\beta-\alpha\beta-\beta}{2}})\xi(r)^{-\frac{1}{2}(n-1)} & \text{if} \hspace{4mm} R^{\beta}\le r\le \frac{3}{2}R^{\beta}; \\ 0 & \text{if}\hspace{4mm} r\ge \frac{3}{2}R^{\beta}.\end{array}\right.\end{equation*}
Our strategy will be similar to that of the proofs of the previous theorems. We will use the expression (\ref{2.4}) again. Fixed $R$, let $K_R$ be the value of the expression (\ref{2.4}) on interval $[0,R]$, then $K_R$ is a constant. Let us analyze the expression (\ref{2.4}) on the interval $[R,R^{\beta}]$. The expression is the same of (\hyperlink{17}{17}), that is 

\vspace{10pt}

\noindent$\displaystyle=\left(\frac{(n\alpha-\alpha-1)^2}{4}-a\alpha(n-1)(n\alpha-2)\right)\int_R^{R^{\beta}}r^{-1}dr$

\noindent$\displaystyle+\frac{1}{2}(n-1)(n\alpha-\alpha-4a\alpha n+4a-1)\int_R^{R^{\beta}}\xi(r)^{-1}\xi'(r)dr$

\noindent$\displaystyle+(n-1)\left(\frac{1}{4}n-an-\frac{1}{4}\right)\int_R^{R^{\beta}}r\xi(r)^{-2}\xi'(r)^2dr$

\noindent$\displaystyle+a\frac{S(F)}{A(F)}\int_{R}^{R^{\beta}}r^{-2\alpha+1}\xi(r)^{-2}dr.\hspace{109.4mm}\hypertarget{22}{(22)}$

\vspace{10pt}

\noindent the relation (\ref{appr}) implies that the second term of (\hyperlink{22}{22}) can be arbitrarily small than the first term of (\hyperlink{22}{22}) from an appropriate choice for $\beta$ large. Since $a>\frac{n-1}{4n}$, the third term of \hyperlink{22}{(22)} is nonpositive and the fourth term is equal to $$a\frac{S(F)}{A(F)}\int_R^{R^{\beta}}r\rho(r)^{-2}dr.$$ Suppose that $$\lim_{r\rightarrow\infty}\int_R^{R^{\beta}}r\rho(r)^{-2}dr=0,$$ then this term is irrelevant in relation to the first, hence, the expression (\hyperlink{22}{22}) can be negatively so large as we want. 

It remains to show that the expression (\ref{2.4}) in the interval $[R^{\beta},\frac{3}{2}R^{\beta}]$ is upper bounded by a constant independent of $\beta$. The expression, by (\hyperlink{2.15}{19}), is upper limited by 

\vspace{10pt}

\noindent$\displaystyle C+(n-1)\left(\frac{1}{4}n-an-\frac{1}{4}\right)\int_{R^{\beta}}^{\frac{3}{2}R^{\beta}}r^{n\alpha-\alpha}f^2\xi(r)^{n-3}\xi'(r)^2dr$

\noindent $\displaystyle+a\frac{S(F)}{A(F)}\int_{R^{\beta}}^{\frac{3}{2}R^{\beta}}r^{n\alpha-3\alpha}\xi(r)^{n-3}f^2dr$

\noindent $\displaystyle-2a\alpha(n-1)\int_{R^{\beta}}^{\frac{3}{2}R^{\beta}}r^{n\alpha-\alpha-1}f^2\xi(r)^{n-2}\xi'(r)dr$

\noindent$\displaystyle+(2-8\alpha)(n-1)R^{-\frac{n\alpha\beta-\alpha\beta+\beta}{2}}\int_{R^{\beta}}^{\frac{3}{2}R^{\beta}}r^{n\alpha-\alpha}f\xi(r)^{\frac{1}{2}(n-1)-1}\xi'(r)dr$,\hspace{55.9mm}\hypertarget{23}{(23)} 

\vspace{10pt}

\noindent where $C$ is a constant independent of $\beta$. Our objective is to show that (\hyperlink{23}{23}) increases arbitrary less than $\log R$. Since $a>\frac{n-1}{4n}$, we have $\frac{1}{4}n-an-\frac{1}{4}<0$. Let $D$ and $H$ be negative fixed constants such that $$D+H=(n-1)\left(\frac{1}{4}n-an-\frac{1}{4}\right),$$ $E=-2a\alpha(n-1)$ and $J=(2-8a)(n-1)$. Since $$\left(\sqrt{-D}r^{\frac{n\alpha-\alpha}{2}}f\xi(r)^{\frac{n-3}{2}}\xi'(r)-\frac{E}{2\sqrt{-D}}r^{\frac{n\alpha-\alpha-2}{2}}f\xi(r)^{\frac{n-1}{2}}\right)^2\ge 0,$$ we have \begin{equation}\label{lt1} Dr^{n\alpha-\alpha}f^2\xi(r)^{n-3}\xi'(r)^2+Er^{n\alpha-\alpha-1}f^2\xi(r)^{n-2}\xi'(r)+\frac{E^2}{4D}r^{n\alpha-\alpha-2}f^2\xi(r)^{n-1}\le 0\tag{24}\end{equation} and since $$\left(\sqrt{-H}r^{\frac{n\alpha-\alpha}{2}}f\xi(r)^{\frac{n-3}{2}}\xi'(r)-\frac{J}{2\sqrt{-H}}R^{-\frac{n\alpha\beta-\alpha\beta-\beta}{2}}r^{\frac{n\alpha-\alpha}{2}}\right)^2\ge 0,$$ we have \begin{equation*} \label{lt2} Hr^{n\alpha-\alpha}f^2\xi(r)^{n-3}\xi'(r)^2+JR^{-\frac{n\alpha\beta-\alpha\beta-\beta}{2}}r^{n\alpha-\alpha}f\xi(r)^{\frac{1}{2}(n-1)-1}\xi'(r)+\frac{J^2}{4H}R^{-n\alpha\beta-\alpha\beta-\beta}r^{n\alpha-\alpha}\le 0\end{equation*}

\vspace{-10pt}

\begin{center} \hypertarget{25}{(25)}\end{center}Integrating (\ref{lt1}) and (\hyperlink{25}{25}) from $R^{\beta}$ to $\frac{3}{2}R^{\beta}$ ans addicting this two inequalities, we have that the expression (\hyperlink{23}{23}) is upper limited by $$C-\frac{E^2}{4D}\int_{R^{\beta}}^{\frac{3}{2}R^{\beta}}r^{n\alpha-\alpha-2}f^2\xi(r)^{n-1}dr-\frac{J^2}{4H}R^{-n\alpha\beta-\alpha\beta-\beta}\int_{R^{\beta}}^{\frac{3}{2}R^{\beta}}r^{n\alpha-\alpha}dr,$$
which is upper limited by a constant because $f(r)\le R^{-\frac{n\alpha\beta-\alpha\beta+\beta}{2}}\xi(r)^{-\frac{1}{2}(n-1)}$ on $[R^{\beta},\frac{3}{2}R^{\beta}]$. Therefore, the combination of the first, third and fourth term of (\hyperlink{23}{23}) is upper limited by a constant. Using that $f(r)\le R^{-\frac{n\alpha\beta-\alpha\beta+\beta}{2}}\xi(r)^{-\frac{1}{2}(n-1)}$ on $[R^{\beta},\frac{3}{2}R^{\beta}]$, we have that $$\int_{R^{\beta}}^{\frac{3}{2}R^{\beta}}r^{n\alpha-3\alpha}\xi(r)^{n-3}f^2dr\le D\int_{R^{\beta}}^{\frac{3}{2}R^{\beta}}r\rho(r)^{-2}dr$$ for a determined constant $D$, where by our previous assumption, the right-hand side of the inequality grows arbitrarily less than $\log r$. Therefore, under this conditions, $M$ can not be $a$-stable for $a>\frac{n-1}{4n}$. Therefore, $$\limsup_{r\rightarrow\infty}(\log R)^{-1}\int_1^Rr\rho(r)^{-2}dr\ge C$$ for some $C>0$.  

If $S(F)\le 0$,  the fourth term of \hypertarget{22}{(22)} and the second term of \hypertarget{23}{(23)} is nonpositive and, assuming that $\alpha>\frac{2}{n}$, to $M$ be $a$-stable, we must to have
\begin{equation*} \begin{split}&\frac{(n\alpha-\alpha-1)^2}{4}-a\alpha(n-1)(n\alpha-2) \ge 0 \\ \Longrightarrow \hspace{2mm}&  \frac{(n\alpha-\alpha-1)^2}{4}-\frac{n-1}{4n}\alpha(n-1)(n\alpha-2) > 0 \\ \Longrightarrow \hspace{2mm} & n(n^2\alpha^2+\alpha^2+1-2n\alpha^2-2n\alpha+2\alpha)-(n\alpha^2-2\alpha)(n^2-2n+1) > 0 \\ \Longrightarrow\hspace{2mm} & -2\alpha n+2\alpha+n>0 \\ \Longrightarrow \hspace{2mm} & \alpha<\frac{n}{2n-2},\end{split}\end{equation*}

\noindent Therefore, $\lim_{r\rightarrow\infty}\rho(r)r^{-\frac{n}{2n-2}}=0$. It finishes the proof.

\noindent $\blacksquare$

\begin{observacao}
    In the case $n=2$, where $F$ has dimension one, we have $\frac{2}{n}=1$ and $n\alpha-2=0$ when $\alpha=1$. The conclusion will be that the volume growth of $M$ satisfies $$\lim_{R\rightarrow\infty}r^{-2-\delta}Vol(B_r(p))=0,$$ for all $\delta>0$ and $p\in M$.
\end{observacao}

\noindent\textbf{Proof of \Cref{th3.1.8}}. Take as test function the same of the \hyperlink{proof3.1.5}{proof of Theorem 3.1.5}, that is, 
$$f(r)=f_{R,\beta}(r)=\begin{cases}\hspace{40mm}0\hspace{40mm} \hspace{4mm} \text{if} \hspace{4mm} 0\le r\le \frac{m-1}{m}R; \\ \hspace{6.5mm}(2R^{-\frac{n\alpha-\alpha+1}{2}}r-R^{-\frac{n\alpha-\alpha-1}{2}})\xi(R)^{-2a(n-1)} \hspace{7.5mm} \hspace{4mm}\text{if} \hspace{4mm} \frac{R}{2}\le r\le R; \\ \hspace{22,9mm}r^{-\frac{n\alpha-\alpha-1}{2}}\xi(r)^{-2a(n-1)}\hspace{26,9mm} \text{if} \hspace{4mm} R\le r\le R^{\beta}; \\ (-2R^{-\frac{n\alpha\beta-\alpha\beta+\beta}{2}}r+3R^{-\frac{n\alpha\beta-\alpha\beta-\beta}{2}})\xi(r)^{-2a(n-1)}\hspace{7mm}\text{if} \hspace{4mm} R^{\beta}\le r\le \frac{3}{2}R^{\beta}; \\ \hspace{40mm}0\hspace{40mm} \hspace{4mm}\text{if}\hspace{4mm} r\ge \frac{3}{2}R^{\beta}.\end{cases}$$
Then the expression (\ref{2.4}) on interval $[R,R^{\beta}]$ is, by (\hyperlink{16}{16}), 

\vspace{10pt}

\noindent$\displaystyle\left(\frac{(n\alpha-\alpha-1)^2}{4}-a\alpha(n-1)(n\alpha-2)\right)\int_R^{R^{\beta}}r^{-1}\xi(r)^{(1-4a)(n-1)}dr$

\noindent$\displaystyle+2a\alpha(n-1)[(1-4a)(n-1)-1]\int_R^{R^{\beta}}\xi(r)^{(1-4a)(n-1)-1}\xi'(r)dr$

\noindent$\displaystyle+a(n-1)[(1-4a)(n-1)-1]\int_R^{R^{\beta}}r\xi(r)^{(1-4a)(n-1)-2}\xi'(r)^2dr$

\noindent$\displaystyle+a\frac{S(F)}{A(F)}\int_{R}^{R^{\beta}}r^{-2\alpha+1}\xi(r)^{(1-4a)(n-1)-2}dr$

\noindent$=\displaystyle \left (\frac{(n\alpha-\alpha-1)^2}{4}-a\alpha(n-1)(n\alpha-2)\right) C_3\int_R^{R^{\beta}}r^{-1}dr$

\noindent$\displaystyle+2a\alpha(n-1)[(1-4a)(n-1)-1]C_4\int_R^{R^{\beta}}\xi(r)^{-1}\xi'(r)dr$

\noindent$\displaystyle+a(n-1)[(1-4a)(n-1)-1]C_5\int_R^{R^{\beta}}r\xi(r)^{-2}\xi'(r)^2dr$

\noindent$\displaystyle+a\frac{S(F)}{A(F)}C_6\int_{R}^{R^{\beta}}r^{-2\alpha+1}dr$, \hspace{114.8mm}\hypertarget{26}{(26)}

\vspace{10pt}

\noindent where $C_3,C_4,C_5$ and $C_6$ are positive constants, that are dependent of $C_1$, $C_2$ and $a$. Note that in (\hyperlink{26}{26}), the first term is arbitrary less than the third term, the second and the fourth term are limited. Therefore, if $(1-4a)(n-1)-1<0\Leftrightarrow a>\frac{n-2}{4(n-1)}$, (\hyperlink{26}{26}) is negative and, in module, so arbitrary large as we want. On the expression (\ref{2.4}) on the interval $[R^{\beta},\frac{3}{2}R^{\beta}]$, the argument is the same of in the \hyperlink{proof3.1.5}{proof of Theorem 3.1.5}, where we can conclude that the expression (\ref{2.4}) on the interval $[R^{\beta},\frac{3}{2}R^{\beta}]$ is upper limited by a constant. Therefore, $L_a$ for $a>\frac{n-2}{4(n-1)}$ is unstable. It finishes the proof.

\noindent\textbf{Proof of \Cref{cit3.12}}. When $\alpha=1$, we have $n\alpha-3\alpha=n\alpha-\alpha-2=n-3$. Redoing all calculations and steps in the \hyperlink{proof3.1.5}{proof of Theorem 3.1.5}, in the analysis in the interval $[R,R^{\beta}]$ taking $\alpha=1$, the expression (\hyperlink{16}{16}) becomes, for $a\ne\frac{1}{4}$:

\noindent$\displaystyle\frac{n-2}{4}(n-2-4a(n-1))\int_R^{R^{\beta}}r^{-1}\xi(r)^{(1-4a)(n-1)}dr$

\noindent$\displaystyle+\frac{2a(n-1)[(1-4a)(n-1)-1]}{(1-4a)(n-1)}[\xi(R^{\beta})^{(1-4a)(n-1)}-\xi(R)^{(1-4a)(n-1)}]$

\noindent$\displaystyle+a(n-1)[(1-4a)(n-1)-1]\int_R^{R^{\beta}}r\xi(r)^{(1-4a)(n-1)-2}\xi'(r)^2dr$

\noindent$\displaystyle+a(n-1)(n-2)\int_R^{R^{\beta}}r^{-1}\xi(r)^{(1-4a)(n-1)-2}dr.\hspace{86.3mm}\hypertarget{27}{(27)}$

\vspace{10pt}

\noindent It is easy to see that, just as in the case $\alpha>1$ (case of \Cref{cit3.8}), the part of integration in the interval $[R^{\beta},\frac{3}{2}R^{\beta}]$ in is limited by a constant that independent of $\beta$, the only difference is that the last term that is $o(1)$ becomes a term $O(1)$ (see the \hyperlink{proof3.1.5}{proof of Theorem 3.1.5}). Returning to the analysis in the interval $[R,R^{\beta}]$, under the hypothesis $0<C_1\le \xi(r)\le C_2<\infty$, we only have to prove that the expression in (\hyperlink{27}{27}) is negative and large in module for $\beta$ large. In (\hyperlink{27}{27}), we have four terms; the second is limited because $C_1\le\xi(r)\le C_2(r)$ and the third is nonpositive for all $a\ge \frac{n-2}{4(n-1)}$. Then, we only have to prove that the combination of the first and fourth terms is negative and, in module, so large as we want. We have:
$$\int_R^{R^{\beta}}r^{-1}\xi(r)^{(1-4a)(n-1)-2}dr\le C_1^{-2}\int_R^{R^{\beta}}r^{-1}\xi(r)^{(1-4a)(n-1)}dr,$$ 
then

\vspace{10pt}

\noindent$\displaystyle\frac{n-2}{4}(n-2-4a(n-1))\int_R^{R^{\beta}}r^{-1}\xi(r)^{(1-4a)(n-1)}dr$

\noindent$\displaystyle+a(n-1)(n-2)\int_R^{R^{\beta}}r^{-1}\xi(r)^{(1-4a)(n-1)-2}dr$

\noindent$\displaystyle \le \left(\frac{(n-2)^2}{4}+a(n-1)(n-2)(-1+C_1^{-2})\right)\int_R^{R^{\beta}}r^{-1}\xi(r)^{(1-4a)(n-1)}dr$.

\vspace{10pt}

\noindent Since $a>\frac{n-2}{4(n-1)(1-C_1^{-2})}$ and $C_1>1$, the last term above is negative. When $\beta\rightarrow\infty$, the integral $\int_R^{R^{\beta}}r^{-1}\xi(r)^{(1-4a)(n-1)}dr$ tends to $\infty$. It ends the proof for $a\ne\frac{1}{4}$. When $a=\frac{1}{4}$, (\hyperlink{16}{16}) becomes

\vspace{10pt}

\noindent$\displaystyle-\frac{n-2}{4}\int_R^{R^{\beta}}r^{-1}dr-\frac{n-1}{2}\int_R^{R^{\beta}}\xi(r)^{-1}\xi'(r)dr-\frac{n-1}{4}\int_R^{R^{\beta}}r\xi(r)^{-2}\xi'(r)^2dr$

\noindent$\displaystyle+\frac{(n-1)(n-2)}{4}\int_R^{R^{\beta}}r^{-1}\xi(r)^{-2}dr$

\noindent $\displaystyle \le \frac{n-2}{4}(-1+(n-1)C_1^{-2})(\beta-1)\log R-\frac{n-1}{2}(\log(\xi(R^{\beta}))-\log(\xi(R)))$,

\vspace{10pt}

\noindent in which will be negative and large in module, if $-1+(n-1)C_1^{-2}<0\Leftrightarrow C_1>(n-1)^{\frac{1}{2}}$. In this conditions, $$\frac{n-2}{4(n-1)(1-C_1^{-2})}<\frac{n-2}{4(n-1)}\frac{1}{1-\frac{1}{n-1}}=\frac{1}{4}.$$

\vspace{10pt}

\noindent Hence, the same conclusions follow for $a=\frac{1}{4}$.

\noindent$\blacksquare$

\noindent\textbf{Proof of \Cref{cit3.13}}.: Similarly to the proof of the previous theorem, we have that $$\int_R^{R^{\beta}}r^{-1}\xi(r)^{(1-4a)(n-1)-2}dr\ge C_2^{-2}\int_R^{R^{\beta}}r^{-1}\xi(r)^{(1-4a)(n-1)}dr$$

\noindent and

\vspace{10pt}

\noindent $\displaystyle\frac{n-2}{4}(n-2-4a(n-1))\int_R^{R^{\beta}}r^{-1}\xi(r)^{(1-4a)(n-1)}dr$

\noindent$\displaystyle+a(n-1)(n-2)\int_R^{R^{\beta}}r^{-1}\xi(r)^{(1-4a)(n-1)-2}dr$

\noindent$\displaystyle \le \left(\frac{(n-2)^2}{4}C_2^2+a(n-1)(n-2)(1-C_2^2)\right)\int_R^{R^{\beta}}r^{-1}\xi(r)^{(1-4a)(n-1)-2}dr$.

\vspace{10pt}

\noindent Since $a<\frac{n-2}{4(n-1)(1-C_2^{-2})}$ and $C_2<1$, the the last term above is negative. When $\beta\rightarrow\infty$, the integral $\int_R^{R^{\beta}}r^{-1}\xi(r)^{(1-4a)(n-1)}dr$ tends to $\infty$. It ends the proof.

\noindent$\blacksquare$

\section{Warped-product minimal hypersurfaces on \texorpdfstring{$\mathbb R^n$}{k}}

There are three typical examples of minimal hypersurfaces in Euclidean spaces that can be written as a warped product. 

The first is any affine hyperplane $H\subset\mathbb R^{n+1}$, which geometrically is the warped product $\mathbb R\times_{\rho}\mathbb S^{n-1}$ where $ \rho(r)=r$, having the topology of $\mathbb R^n$, being stable as a minimal hypersurface of $\mathbb R^{n+1}$ and, more than it, has the property of minimizing the area for any variation of compact support. When $n\le 5$, the affine hyperplane is the only minimal hypersurface which is the graph of a smooth function $f:\mathbb R^n\rightarrow\mathbb R$ (the problem of determining whether the graph of a function of $\mathbb R^n$ is a minimal hypersurface on $\mathbb R^{n+1}$ is known as the \textit{Bernstein problem}. For more details, see \autocite{bernstein1917theoreme}, \autocite{almgren1966some} and \autocite{chodosh2021stablev}). According to our notation, because $H$ is a flat manifold, it is $a$-stable for all $a\in\mathbb R$.

The second example is the $n$-dimensional catenoid in $\mathbb R^{n+1}$. In dimension two, the catenoid can be parameterized by $x=c\cosh(v/c)\cos u$, $y=c\cosh(v/c)\sin(u)$, $z=v$, where $c\ne 0$, $v\in\mathbb R$ and $y\in [0,2\pi)$ and is characterized by being a minimal surface in $\mathbb R^3$ of revolution (different from plane) and index one (see \autocite{tam2009stability}), which can be written as a warped product $\mathbb R\times_{\rho}\mathbb S^1$. In higher dimensions, the catenoid is a minimal hypersurface of $\mathbb R^{n+1}$ that inherits certain properties from the two-dimensional catenoid. For example, it is diffeomorphic to the cylinder $\mathbb R\times\mathbb S ^{n-1}$ and has index one. Perhaps the main difference for dimension two is that in higher dimensions, the catenoid is limited in one of the directions of $\mathbb R^{n+1}$. The $n$-dimensional catenoid is $\frac{n-2}{n}$-stable (see \autocite{tam2009stability}).

The third example is the minimal cones in $\mathbb R^{n+1}$, characterized as being a warped product $[0,\infty)\times_{\rho}F$, where $\rho(r)=r$ and $F$ is a minimal submanifold of $\mathbb S^n$. A particular case is the Simons cones, introduced by J. Simons on \autocite{simons1968minimal}, and also worked by E. Bombieri, E. De Giorgi, and E. Giusti on \autocite{bombieri1969minimal}. They are characterized as the singular hypersurface of $\mathbb R^{2m}$: $x_1^2+\dots+x_m^2=x_{m+1}^2+\dots+x_{2m}^2$. These cones can be written as the warped product $[0,\infty)\times_{\rho}\left(\frac{1}{\sqrt{2}}\mathbb S^{n-1}\times \frac{1}{\sqrt {2}}\mathbb S^{n-1}\right)$, where $\rho(r)=r$ and is a minimal hypersurface of $\mathbb R^{2m}$, more precisely, these cones are $\frac{(2m-3)^2 }{8(m-1)}$-stable according to our definition of $a$-stability. Note that $\frac{(2n-3)^2}{8(n-1)}>1$, if $n\ge 4$. In particular, a Simons cone is a minimal stable hypersurface of $\mathbb R^{2n}$ if and only if $2n\ge 8$.  

\begin{proposicao}\label{propcone}An $a$-stable minimal cone $C^n$ of $\mathbb R^{n+1}$, with $a>\max\left\{1,\frac{(n-2)^2}{4(n-1)}\right\}$, is flat.\end{proposicao} 

\noindent Simons in \autocite{simons1968minimal} has shown that a minimal $1$-stable cone $C^n\subset\mathbb R^{n+1}$ is flat if $n\le 6$. The \Cref{propcone} and the previous examples of the Simons cone shows that, for $n\ge 7$, the value $\frac{(n-2)^2}{4(n-1)}$ 
 is the greatest possible value for $a$ such that a $n$-dimensional minimal cone in $\mathbb R^{n+1}$ can be $a$-stable.

\begin{proposicao}{(Simons \autocite{simons1968minimal})} Suppose that $M\rightarrow\mathbb R^n$ is a minimal immersion. Then, the second fundamental for $A$ satisfies $$\frac{1}{2}\Delta|A|^2+|A|^4=|\nabla A|^2$$

\noindent on $M$.
\end{proposicao}

For a demonstration, see (\autocite{chodosh2021stablev}, Proposition 8.13). It implies \begin{equation}\label{2}|A|\Delta |A|+|A|^4=|\nabla A|^2-|\nabla|A||^2.\tag{28}\end{equation}

\noindent \textbf{Proof of \Cref{propcone}}. Let $A_F$ be the second fundamental form of the immersion $F\rightarrow \mathbb S^n$ and $A$ be the second fundamental form of the immersion $C\rightarrow\mathbb R^{n+1}$. The cone $C$ can be seen as the set $$C=rF=\{rx,x\in F, r\ge 0\}.$$

\noindent Let $rx\in C$, $r\ge 0$, $x\in F$, by a direct calculus, $A(\partial r,\cdot)=0$ and for $X,Y\in T_xF$, $A(X,Y)=r^{-1}A_F(X,Y)$. Then \begin{equation*}\begin{split}(\nabla_{\partial r}A)(X,Y)=&\nabla_{\partial r}(A(X,Y))-A(\nabla_{\partial r}X,Y)-A(X,\nabla_{\partial r}Y) \\ = &-r^{-2}A_F(X,Y) \\ = & -r^{-1}A(X,Y),\end{split}\end{equation*}

\noindent where $\nabla$ is the Levi-Civita connection of $C$. We can then consider $A=r^{-1}A_F$ and $\nabla_{\partial r}A=-r^{-1}A$. Being $\{e_1,\dots,e_{n-1}\}$ an orthonormal frame of $F$ and $e_n=\partial r$ such that $\{e_1,\dots,e_n\}$ on a given point $p=(r,x)$ is geodesic and diagonalizing $A$, then $$|\nabla A|^2=\sum_{i,j,k=1}^n\left[\left(\nabla_{e_i}A\right)\left(e_j,e_k\right)\right]^2$$ and 
\begin{equation*}\begin{split}|\nabla |A||^2 & =\sum_{i=1}^n\langle \nabla|A|,e_i\rangle^2 \\ & =\sum_{i=1}^n\left\langle \nabla\left(\sum_{j,k=1}^nA\left(e_j,e_k\right)^2\right)^{\frac{1}{2}},e_i\right\rangle^2 \\ &=|A|^{-2}\sum_{i,j,k=1}^nA(e_j,e_k)^2[e_i(A(e_j,e_k))]^2 \\ &=|A|^{-2}\sum_{i,j=1}^nA(e_j,e_j)^2[e_i(A(e_j,e_j))]^2 \\ &=|A|^{-2}\sum_{i,j-1}^nA(e_j,e_j)^2[(\nabla_{e_i}A)(e_j,e_j)]^2 \\ & \le |A|^{-2}\sum_{i=1}^n\left(\sum_{j=1}^nA(e_j,e_j)^2\right)\left(\sum_{j=1}^n[(\nabla_{e_i}A)(e_j,e_j)]^2\right)\\ &=\sum_{i,j=1}^n[\nabla_{e_i}A(e_j,e_j)]^2.\end{split}\end{equation*}

\noindent Hence \begin{equation*}\begin{split}|\nabla A|^2-|\nabla|A||^2 & \ge 2\sum_{i=1}^n\sum_{1\le j<k\le n}[(\nabla_{e_i}A)(e_j,e_k)]^2 \\ & \ge 2\sum_{i=1}^n\sum_{j=1}^{n-1}[(\nabla e_i A)(e_j,e_n)]^2.\end{split}\end{equation*}

\noindent By Codazzi equation, $(\nabla_{e_i}A)(e_j,e_n)=(\nabla_{e_n}A)(e_i,e_j)$. Therefore, \begin{equation*}\begin{split}|\nabla A|^2-|\nabla|A||^2 & \ge 2\sum_{i=1}^n\sum_{j=1}^{n-1}[(\nabla e_n A)(e_i,e_j)]^2 \\ & =2r^{-2}|A|^2.\end{split}\end{equation*}

\noindent It and (\ref{2}) imply $$|A|\Delta|A|+|A|^4\ge 2r^{-2}|A|^2.$$

\noindent For $\varphi\in C_c^{\infty}(C\setminus\{0\})$, multiplying the inequality above by $|A|^{\frac{1-a}{a}}\varphi^2$ and integrating on $C$, we have
\begin{equation}\label{3}\begin{split} 2\int_C|A|^{\frac{1+a}{a}}r^{-2}\varphi^2\le &  \int_C|A|^{\frac{1}{a}}\varphi^2\Delta |A|+|A|^{\frac{1+3a}{a}}\varphi^2\\ = & \int_C|A|^{\frac{1+3a}{a}}\varphi^2-\frac{1}{a}|A|^{\frac{1-a}{a}}\varphi^2|\nabla|A||^2-2|A|^{\frac{1}{a}}\varphi\langle\nabla\varphi,\nabla|A|\rangle.\end{split}\tag{29}\end{equation}

\noindent Substituting $|A|^{\frac{a+1}{2a}}\varphi$ on $a$-stability inequality (recall that $|A|^2=-S$ on minimal hypersurfaces of $\mathbb R^n$), we have:
\begin{equation*}\begin{split}a\int_C |A|^{\frac{1+3a}{a}}\varphi^2\le & \int_C|\nabla (|A|^{\frac{a+1}{2a}}\varphi)|^2 \\ = & \int_C\left(\frac{a+1}{2a}\right)^2|A|^{\frac{1-a}{a}}\varphi^2|\nabla |A||^2+|A|^{\frac{a+1}{a}}|\nabla\varphi|^2 +\frac{a+1}{a}|A|^{\frac{1}{a}}\varphi\langle\nabla|A|,\nabla\varphi\rangle\end{split}\end{equation*}

\begin{center}\hypertarget{3.15}{(30)}\end{center}
 Multiplying (\ref{3}) by $a$ and adding to (\hyperlink{3.15}{30}), we find
\begin{equation*}\begin{split}\displaystyle 2a\int_C|A|^{\frac{a+1}{a}}r^{-2}\varphi^2\le &  \int_C\left[\left(\frac{a+1}{2a}\right)^2-1\right]|A|^{\frac{1-a}{a}}\varphi^2|\nabla|A||^2 \\ &+\frac{-2a^2+a+1}{a}|A|^{\frac{1}{a}}\varphi\langle \nabla |A|,\nabla\varphi\rangle+|A|^{\frac{a+1}{a}}|\nabla\varphi|^2 \\ = & \int_C\frac{(-3a-1)(a-1)}{4a^2}|A|^{\frac{1-a}{a}}\varphi^2|\nabla|A||^2 \\ &+\frac{(-2a-1)(a-1)}{a}|A|^{\frac{1}{a}}\varphi\langle \nabla |A|,\nabla\varphi\rangle+|A|^{\frac{a+1}{a}}|\nabla\varphi|^2.\hspace{1cm}\hypertarget{3.12}{(31)}\end{split}\end{equation*}

\noindent Take $\varphi=\varphi(r)$ only dependent of $r$. Since $|A|=r^{-1}|A_F|$, then $$\nabla\varphi=\varphi'(r)\partial_r, \hspace{3mm}\langle \nabla |A|,\nabla\varphi\rangle=-r^{-2}\phi'(r)|A_F|\hspace{3mm} \text{and}\hspace{3mm} |\nabla|A||^2\ge r^{-4}|A_F|^2.$$ Hence (\hyperlink{3.12}{31}) with the hypothesis $a>1$ implies
\begin{equation*}\begin{split} 2a\int_C|A_F|^{\frac{a+1}{a}}r^{-\frac{3a+1}{a}}\varphi^2\le & \int_C-\frac{3a^2-2a-1}{4a^2}r^{-\frac{3a+1}{a}}|A_F|^{\frac{a+1}{a}}\varphi^2\\ &+\frac{2a^2-a-1}{a}|A_F|^{\frac{a+1}{a}}r^{-\frac{2a+1}{a}}\varphi\varphi'+r^{-\frac{a+1}{a}}|A_F|^{\frac{a+1}{a}}\varphi'^2.\end{split}\end{equation*}
\noindent Since $dV_C=r^{n-1}dV_Fdr$, the last inequality is equivalent to:
$$0\le \int_F|A_F|^{\frac{a+1}{a}}\int_0^{\infty}-\frac{8a^3+3a^2-2a-1}{4a^2}r^{\frac{na-4a-1}{a}}\varphi^2+\frac{2a^2-a-1}{a}r^{\frac{na-3a-1}{a}}\varphi\varphi'+r^{\frac{na-2a-1}{a}}\varphi'^2.$$

\noindent If $\int_F|A_F|^{\frac{a+1}{a}}=0$, then $|A_F|\equiv 0$ and $|A|\equiv 0$, therefore, $C$ is flat. Thus, suppose that $C$ is non flat, then
\begin{equation}\label{5}0\le \int_0^{\infty}-\frac{8a^3+3a^2-2a-1}{4a^2}r^{\frac{na-4a-1}{a}}\varphi^2+\frac{2a^2-a-1}{a}r^{\frac{na-3a-1}{a}}\varphi\varphi'+r^{\frac{na-2a-1}{a}}\varphi'^2.\tag{32}\end{equation}

\noindent By an argument of approximation, we can take $\varphi\in C_c^{0,1}(C)$ defined by: $$\varphi(r)=\varphi_R(r)=\begin{cases} \hspace{22.3mm}0\hspace{22.3mm}\hspace{4mm}\text{if}\hspace{2mm} 0\le r\le\frac{1}{2};\\ \hspace{21.2mm}2r\hspace{21.2mm}\hspace{4mm} \text{if}\hspace{2mm} \frac{1}{2}\le r\le 1; \\ \hspace{15.3mm}r^{\frac{-na+3a+1}{2a}}\hspace{15.3mm}\hspace{4mm} \text{if} \hspace{2mm}1\le r\le R; \\ \hspace{0.5mm}-R^{\frac{-na+a+1}{2a}}r+2R^{\frac{-na+3a+1}{2a}}\hspace{0.5mm}\hspace{4mm}\text{if} \hspace{2mm}R\le r\le 2R;\\ \hspace{22.3mm}0\hspace{22.3mm}\hspace{4mm}\text{if}\hspace{2mm} r>2R.\end{cases}$$

\noindent The integral (\ref{5}) on the interval $[0,1]$ is a constant independent of $R$. The integral on the interval $[1,R]$ becomes 

\vspace{4mm}

\noindent$\displaystyle\int_1^{R}\left[\frac{-8a^3-3a^2+2a+1}{4a^2}+\frac{2a^2-a-1}{a}\frac{-na+3a+1}{2a}+\frac{(na-3a-1)^2}{4a^2}\right]r^{-1}dr$

\noindent$\displaystyle=\left(\frac{-8a^3-3a^2+2a+1}{4a^2}+\frac{-4na^3+12a^3+4a^2+2na^2-6a^2-2a+2na-6a-2}{4a^2}\right.$

$\displaystyle\left.+\frac{n^2a^2+9a^2+1-6na^2-2na+6a}{4a^2}\right)\log R$

\noindent $=\displaystyle \frac{-4na^3+4a^3+n^2a^2-4na^2+4a^2}{4a^2}\log R$

\noindent$\displaystyle =\left(-na+a+\frac{n^2}{4}-n+1\right)\log R$.

\vspace{4mm}

\noindent Hence, the integral of the expression (\ref{5}) on interval $[1,R]$ will be negative and a multiple of $\log R$, if
$$-na+a+\frac{n^2}{4}-n+1<0\Longleftrightarrow a(1-n)<-\frac{(n-2)^2}{4}\Longleftrightarrow a>\frac{(n-2)^2}{4(n-1)}.$$
\noindent Let us check that the value of expression (\ref{5}) on the interval $[R,2R]$ is constant and independent of $R$. Note that:

\vspace{4mm}

\noindent$\displaystyle\int_R^{2R}r^{\frac{na-4a-1}{a}}\varphi(r)^2dr=R^{\frac{-na+a+1}{a}}\int_R^{2R}\left(r^{\frac{na-2a-1}{a}}-4Rr^{\frac{na-3a-1}{a}}+4R^2r^{\frac{na-4a-1}{a}}\right)dr=C_1$;

\noindent$\displaystyle\int_R^{2R}r^{\frac{na-3a-1}{a}}\varphi(r)\varphi'(r)dr=R^{\frac{-na+a+1}{a}}\int_R^{2R}\left(-r^{\frac{na-2a-1}{a}}+2Rr^{\frac{na-3a-1}{a}}\right)dr=C_2$;

\noindent$\displaystyle\int_R^{2R}r^{\frac{na-2a-1}{a}}\varphi'(r)^2dr=R^{\frac{-na+a+1}{a}}\int_R^{2R}r^{\frac{na-2a-1}{a}}dr=C_3$,

\vspace{4mm}

\noindent where $C_1$, $C_2$ and $C_3$ are independent of $R$. Hence, the expression of (\ref{5}) on the interval $[R,2R]$ is constant and independent of $R$.

Therefore, for $R$ large, $\varphi_R$ negatives (\ref{5}) and $C$ cannot be $a$-stable for $a>\frac{(n-2)^2}{4(n-1)}$, if $C$ is non flat.

\noindent $\blacksquare$

\nocite{*}

\printbibliography

@article{kawai1988operator,
  title={Operator $\Delta $-aK on surfaces},
  author={Kawai, Shigeo},
  journal={Hokkaido Mathematical Journal},
  volume={17},
  number={2},
  pages={147--150},
  year={1988},
  publisher={Hokkaido University, Department of Mathematics}
}

@article{fischer1980structure,
  title={The structure of complete stable minimal surfaces in 3-manifolds of non-negative scalar curvature},
  author={Fischer-Colbrie, Doris and Schoen, Richard},
  journal={Communications on Pure and Applied Mathematics},
  volume={33},
  number={2},
  pages={199--211},
  year={1980},
  publisher={Citeseer}
}

@article{do1979stable,
  title={Stable complete minimal surfaces in $\mathbb R^3$ are planes},
  author={do Carmo, M and Peng, CK},
  journal={Bulletin of the American Mathematical Society},
  volume={1},
  number={6},
  pages={903--906},
  year={1979}
}

@article{fischer1985complete,
  title={On complete minimal surfaces with finite Morse index in three manifolds},
  author={Fischer-Colbrie, Doris},
  journal={Inventiones mathematicae},
  volume={82},
  number={1},
  pages={121--132},
  year={1985},
  publisher={Springer}
}

@article{nayatani1990morse,
  title={On the Morse index of complete minimal surfaces in Euclidean},
  author={Nayatani, Shin},
  journal={Osaka Journal of Mathematics},
  volume={27},
  number={2},
  pages={441--451},
  year={1990}
}

@article{morabito2009index,
  title={Index and nullity of the Gauss map of the Costa-Hoffman-Meeks surfaces},
  author={Morabito, Filippo},
  journal={Indiana University mathematics journal},
  pages={677--707},
  year={2009},
  publisher={JSTOR}
}

@article{lopez1991embedded,
  title={On embedded complete minimal surfaces of genus zero},
  author={L{\'o}pez, Francisco J and Ros, Antonio},
  journal={Journal of Differential Geometry},
  volume={33},
  number={1},
  pages={293--300},
  year={1991},
  publisher={Lehigh University}
}

@article{costa1989uniqueness,
  title={Uniqueness of minimal surfaces embedded in R3 with total curvature 12$\pi$},
  author={Costa, C},
  journal={J. Differential Geom},
  volume={30},
  number={3},
  pages={597--618},
  year={1989}
}

@article{chodosh2018topology,
  title={On the topology and index of minimal surfaces II},
  author={Chodosh, Otis and Maximo, Davi},
  journal={Journal of Differential Geometry},
  volume={123},
  number={3},
  pages={431--459},
  year={2023},
  publisher={Lehigh University}
}

@book{osserman2013survey,
  title={A survey of minimal surfaces},
  author={Osserman, Robert},
  year={2013},
  publisher={Courier Corporation}
}

@article{choe1990index,
  title={Index, vision number and stability of complete minimal surfaces},
  author={Choe, Jaigyoung},
  journal={Archive for Rational Mechanics and Analysis},
  volume={109},
  pages={195--212},
  year={1990},
  publisher={Springer-Verlag}
}

@article{ritore1997index,
  title={Index one minimal surfaces in flat three space forms},
  author={Ritor{\'e}, Manuel},
  journal={Indiana University mathematics journal},
  pages={1137--1153},
  year={1997},
  publisher={JSTOR}
}

@book{chen2017differential,
  title={Differential geometry of warped product manifolds and submanifolds},
  author={Chen, Bang-Yen},
  year={2017},
  publisher={World Scientific}
}

@article{perelman2002entropy,
  title={The entropy formula for the Ricci flow and its geometric applications},
  author={Perelman, Grisha},
  journal={arXiv preprint math/0211159},
  year={2002}
}

@article{perelman2003ricci,
  title={Ricci flow with surgery on three-manifolds},
  author={Perelman, Grisha},
  journal={arXiv pre-print math/0303109},
  year={2003}
}

@article{berard2014inverse,
  title={Inverse spectral positivity for surfaces},
  author={B{\'e}rard, Pierre and Castillon, Philippe},
  journal={Revista matem{\'a}tica iberoamericana},
  volume={30},
  number={4},
  pages={1237--1264},
  year={2014}
}

@article{barbosa1976size,
  title={On the size of a stable minimal surface in R3},
  author={Barbosa, Joao Lucas and Carmo, M do},
  journal={American Journal of Mathematics},
  pages={515--528},
  year={1976},
  publisher={JSTOR}
}

@article{peetre1957generalization,
  title={A generalization of Courant's nodal domain theorem},
  author={Peetre, Jaak},
  journal={Mathematica Scandinavica},
  pages={15--20},
  year={1957},
  publisher={JSTOR}
}

@article{nayatani1993morse,
  title={Morse index and Gauss maps of complete minimal surfaces in Euclidean 3-space},
  author={Nayatani, Shin},
  journal={Commentarii Mathematici Helvetici},
  volume={68},
  number={1},
  pages={511--537},
  year={1993},
  publisher={Springer}
}

@inproceedings{pogorelov1981stability,
  title={On the stability of minimal surfaces},
  author={Pogorelov, Aleksei Vasil'evich},
  booktitle={Doklady Akademii Nauk},
  volume={260},
  number={2},
  pages={293--295},
  year={1981},
  organization={Russian Academy of Sciences}
}

@inproceedings{montiel2006schrodinger,
  title={Schr{\"o}dinger operators associated to a holomorphic map},
  author={Montiel, Sebasti{\'a}n and Ros, Antonio},
  booktitle={Global Differential Geometry and Global Analysis: Proceedings of a Conference held in Berlin, 15--20 June, 1990},
  pages={147--174},
  year={2006},
  organization={Springer}
}

@incollection{nayatani1992morse,
  title={Morse index of complete minimal surfaces},
  author={Nayatani, Shin},
  booktitle={The Problem Of Plateau: A Tribute to Jesse Douglas and Tibor Rad{\'o}},
  pages={181--189},
  year={1992},
  publisher={World Scientific}
}

@article{nayatani1990lower,
  title={Lower bounds for the Morse index of complete minimal surfaces in},
  author={Nayatani, Shin},
  journal={Osaka Journal of Mathematics},
  volume={27},
  number={2},
  pages={453--464},
  year={1990}
}

@article{chodosh2021stablev,
  title={Stable minimal surfaces and positive scalar curvature lecture notes for math 258, Stanford, Fall 2021}, author={Chodosh, Otis},
  journal={},
  year={2021}
}

@article{chodosh2021stable,
  title={Stable minimal hypersurfaces in $\mathbb R^4$},
  author={Chodosh, Otis and Li, Chao},
  journal={arXiv preprint arXiv:2108.11462},
  year={2021}
}

@article{dobarro2005curvature,
  title={Curvature of multiply warped products},
  author={Dobarro, Fernando and {\"U}nal, B{\"u}lent},
  journal={Journal of Geometry and Physics},
  volume={55},
  number={1},
  pages={75--106},
  year={2005},
  publisher={Elsevier}
}

@article{li2023metrics,
  title={Metrics with $\lambda_1 (-\Delta+ k R)\ge 0$ and Flexibility in the Riemannian Penrose Inequality},
  author={Li, Chao and Mantoulidis, Christos},
  journal={Communications in Mathematical Physics},
  volume={401},
  number={2},
  pages={1831--1877},
  year={2023},
  publisher={Springer}
}

@article{meeks2006liouville,
  title={Liouville-type properties for embedded minimal surfaces},
  author={Meeks, William H and P{\'e}rez, Joaqu{\'\i}n and Ros, Antonio},
  journal={Communications in Analysis and Geometry},
  volume={14},
  number={4},
  pages={703--723},
  year={2006},
  publisher={International Press of Boston}
}

@article{tam2009stability,
  title={Stability properties for the higher dimensional catenoid in $\mathbb R^{n+1}$},
  author={Tam, Luen-Fai and Zhou, Detang},
  journal={Proceedings of the American Mathematical Society},
  volume={137},
  number={10},
  pages={3451--3461},
  year={2009}
}

@article{da1987stability,
  title={Stability of complete noncompact surfaces with constant mean curvature},
  author={Da Silveira, Alexandre M},
  journal={Mathematische Annalen},
  volume={277},
  pages={629--638},
  year={1987},
  publisher={Springer}
}

@article{lopez1989complete,
  title={Complete minimal surfaces with index one and stable constant mean curvature surfaces},
  author={L{\'o}pez, Francisco J and Ros, Antonio},
  journal={Comment. Math. Helv},
  volume={64},
  number={1},
  pages={34--43},
  year={1989}
}

@article{espinar2011colding,
  title={A Colding-Minicozzi stability inequality and its applications},
  author={Espinar, Jos{\'e} and Rosenberg, Harold},
  journal={Transactions of the American Mathematical Society},
  volume={363},
  number={5},
  pages={2447--2465},
  year={2011}
}

@article{colding2002estimates,
  title={Estimates for parametric elliptic integrands},
  author={Colding, Tobias H and Minicozzi, William P},
  journal={International Mathematics Research Notices},
  volume={2002},
  number={6},
  pages={291--297},
  year={2002},
  publisher={OUP}
}

@article{chodosh2016topology,
  title={On the topology and index of minimal surfaces},
  author={Chodosh, Otis and Maximo, Davi},
  journal={Journal of Differential Geometry},
  volume={104},
  number={3},
  pages={399--418},
  year={2016},
  publisher={Lehigh University}
}

@article{bombieri1969minimal,
  title={Minimal cones and the Bernstein problem},
  author={Bombieri, E and De Giorgi, E and Giusti, E},
  journal={Ennio De Giorgi},
  volume={291},
  year={1969},
  publisher={Springer}
}

@article{simons1968minimal,
  title={Minimal varieties in Riemannian manifolds},
  author={Simons, James},
  journal={Annals of Mathematics},
  pages={62--105},
  year={1968},
  publisher={JSTOR}
}

@article{cao1997structure,
  title={The structure of stable minimal hypersurfaces in $\mathbb R^{n+1}$},
  author={Cao, Huai-Dong and Shen, Ying and Zhu, Shunhui},
  journal={Mathematical Research Letters},
  volume={4},
  number={5},
  pages={637--644},
  year={1997},
  publisher={International Press of Boston}
}

@article{yau1976harmonic,
  title={Harmonic Maps and the Topology of Stable Hypersurfaces and Manifolds with Non-negative Ricci Curvature},
  author={Yau, Shing Tung and Schoen, Richard},
  journal={Commentarii mathematici Helvetici},
  volume={51},
  pages={333--342},
  year={1976}
}

@article{bernstein1917theoreme,
  title={Sur un th{\'e}or{\`e}me de g{\'e}om{\'e}trie et ses applications aux {\'e}quations aux d{\'e}riv{\'e}es partielles du type elliptique},
  author={Bernstein, Serge},
  journal={Comm. de la Soc. Math. de Kharkov (2{\'e}me s{\'e}r.)},
  volume={15},
  pages={38--45},
  year={1917}
}

@article{almgren1966some,
  title={Some interior regularity theorems for minimal surfaces and an extension of Bernstein's theorem},
  author={Almgren, Frederick J},
  journal={Annals of Mathematics},
  pages={277--292},
  year={1966},
  publisher={JSTOR}
}

@article{bishop1969manifolds,
  title={Manifolds of negative curvature},
  author={Bishop, Richard L and O’Neill, Barrett},
  journal={Transactions of the American Mathematical Society},
  volume={145},
  pages={1--49},
  year={1969}
}

@article{neto2014topicos,
  title={T{\'o}picos de Geometria Diferencial},
  author={Neto, Antonio Caminha Muniz},
  journal={Sociedade Brasileira de Matem{\'a}tica},
  year={2014}
}

\end{document}